\documentclass[aap,seceqn,nameyear,MSNbibl,dvips]{arximspdf}

\doi{10.1214/10-AAP714}
\volume{21}
\issue{5}
\pubyear{2011}
\firstpage{1663}
\lastpage{1693}

\makeatletter

\newtheorem{theorem}{Theorem}
\newproclaim{definition}{Definition}
\newtheorem{lemma}{Lemma}[section]
\newtheorem{proposition}{Proposition}
\newproclaim{remark}{Remark}[section]

\newcommand{\varlimsup}{\mathop{\overline{\lim}}}
\newcommand{\varliminf}{\mathop{\underline{\lim}}}

\newcommand{\vep}{{\varepsilon}}

\newcommand{\cB}{{\mathcal B}}
\newcommand{\cL}{{\mathcal L}}
\newcommand{\BbbR}{{\mathbb R}}
\newcommand{\BbbN}{{\mathbb N}}

\newcommand{\LM}{{\mathbb A}}
\newcommand{\bB}{{\mathbf B}}

\newcommand{\BbbP}{{\mathbb P}}
\newcommand{\cF}{{\mathcal F}}
\newcommand{\cP}{{\mathcal P}}

\makeatother

\begin{document}
\begin{frontmatter}

\title{Existence and uniqueness of solutions to the inverse boundary
crossing problem for diffusions}
\runtitle{Inverse boundary crossing problem}

\begin{aug}
\author[A]{\fnms{Xinfu} \snm{Chen}\thanksref{t1}\ead[label=e1]{xinfu@pitt.edu}},
\author[B]{\fnms{Lan} \snm{Cheng}\ead[label=e2]{Lan.Cheng@fredonia.edu}},
\author[A]{\fnms{John} \snm{Chadam}\thanksref{t3}\ead[label=e3]{chadam@pitt.edu}} and
\author[C]{\fnms{David} \snm{Saunders}\corref{}\thanksref{t4}\ead[label=e4]{dsaunders@uwaterloo.ca}}
\runauthor{Chen, Cheng, Chadam and Saunders}
\affiliation{University of Pittsburgh, SUNY, Fredonia, University of
Pittsburgh and University of Waterloo}
\address[A]{X. Chen\\
J. Chadam\\
Department of Mathematics\\
University of Pittsburgh\\
Pittsburgh, Pennsylvania 15260\\
USA\\
\printead{e1}\\
\phantom{E-mail: }\printead*{e3}}
\address[B]{L. Cheng\\
Department of Mathematical Sciences\\
SUNY Fredonia\\
Fredonia New York 14063\\
USA\\
\printead{e2}}
\address[C]{D. Saunders\\
Department of Statistics and Actuarial Science\\
University of Waterloo\\
Waterloo, Ontario N2L 3G1\\
Canada\\
\printead{e4}}
\end{aug}

\thankstext{t1}{Supported by NSF Grant DMS-05-04691.}
\thankstext{t3}{Supported by NSF Grant DMS-07-07953.}
\thankstext{t4}{Supported by an NSERC Discovery grant.}

\received{\smonth{4} \syear{2009}}
\revised{\smonth{4} \syear{2010}}

%
\begin{abstract}
We study the inverse boundary crossing problem for diffusions. Given a diffusion
process $X_{t}$, and a survival distribution $p$ on $[0,\infty)$, we
demonstrate
that there exists a boundary $b(t)$ such that $p(t) = \mathbb{P}[\tau
> t]$,
where $\tau$ is the first hitting time of $X_{t}$ to the boundary
$b(t)$. The approach
taken is analytic, based on solving a parabolic variational inequality
to find $b$.
Existence and uniqueness of the solution to this variational inequality
were proven
in earlier work. In this paper, we demonstrate that the resulting
boundary $b$ does
indeed have $p$ as its boundary crossing distribution. Since little is
known regarding
the regularity of $b$ arising from the variational inequality, this
requires a detailed
study of the problem of computing the boundary crossing distribution of
$X_{t}$ to a
rough boundary. Results regarding the formulation of this problem in
terms of weak solutions
to the corresponding Kolmogorov forward equation are
presented.
\end{abstract}

%
\begin{keyword}[class=AMS]
\kwd{60J60}
\kwd{35R35}.
\end{keyword}
\begin{keyword}
\kwd{Inverse boundary crossing problem}
\kwd{Brownian motion}
\kwd{diffusion processes}
\kwd{viscosity solutions}
\kwd{variational inequalities}.
\end{keyword}

\end{frontmatter}

\section{Introduction}\label{intro}
Let $\{B_t\}_{t\geq0}$ be a standard Brownian motion defined
on a filtered probability space $(\Omega,\BbbP,\{\cF_t\}_{t\geq0})$
satisfying
the usual conditions.
We consider a diffusion process $\{X_t\}_{t\geq0}$ defined by the stochastic
differential equation
\[
dX_t = \mu(X_t,t) \,dt + \sigma(X_t,t)
\,dB_t\qquad
\forall t>0,
\]
where
$\mu,\sigma\dvtx\mathbb{R}\times\mathbb{R}_{+}\to\mathbb{R}$
are smooth bounded functions, with bounded
derivatives\setcounter{footnote}{3}\footnote{The main results in the paper are true
under less restrictive assumptions on the coefficients, although we
have not sought to
determine the minimal assumptions under which they hold.
The stated assumptions ensure that the applications of various results from
the theory of PDEs, a change of variables, and an application of It{\^
{o}}'s lemma are all valid.}
and $\inf_{\BbbR\times\BbbR_+}\sigma>0$.
We assume that $X_{0}$ is independent of $B$, and has initial distribution
$\mathbb{P}(X_{0} \leq x) = p_{0}(x,0)$, with density $\rho_{0}(x,0)$.
For $y\in\BbbR$ and $t\geq s$, we further denote by
$F(y,s;\cdot,t)$ and $\rho(y,s;\cdot,t)$ the transition distribution and
density, respectively, of $X_t$ given $X_s=y$:
\[
F(y,s; x,t):= \BbbP( X_t\leq x | X_s=y),\qquad
\rho(y,s;x,t):=\frac{\partial F(y,s;x,t)}{\partial x}.
\]
In the sequel, we denote by $\rho_0(\cdot,t)$ and $p_0(\cdot,t)$
the density\footnote{At $t=0$, $\rho_{0}$ is interpreted as the
distributional derivative of $p_{0}$.}
and cumulative distribution of $X_t$:
\begin{eqnarray*}
p_0(x,t) &:=& \BbbP(X_t\leq x) =\int_\BbbR
\rho_0(y,0)F(y,0;x,t)\,dy,\\
\rho_0(x,t)&:=&\frac{\partial p_0(x,t)}{\partial x}
:=\int_\BbbR\rho_0(y,0) \rho(y,0;x,t) \,dy.
\end{eqnarray*}

For a
given function $b\dvtx\mathbb{R}_{+} \to[-\infty,\infty)$, the first
boundary crossing time
$\hat\tau$ is defined to be
%
%
\begin{equation}\label{BCDefinition}
\hat\tau= \inf\{ t > 0 | X_{t} \leq b(t)\}.
\end{equation}
We shall also have occasion to consider the related, but less commonly
used, first
time that $X_{t}$ goes strictly below $b$:
%
%
\begin{equation} \label{TauDefinition}
\tau= \inf\{ t > 0 | X_{t} < b(t) \}.
\end{equation}

We are interested in the following two problems.
\begin{enumerate}
\item\textit{The boundary crossing problem}: for a given function
$b\dvtx\mathbb{R}_{+}\to[-\infty,\infty)$, compute the survival
distribution of the first
time that $X$ crosses $b$; that is,
%
%
\begin{equation} \label{pDef}
p(t)=\mathbb{P}(\hat\tau\geq t).
\end{equation}
In this case, we denote $p = \cP[b]$.
\item\textit{The inverse boundary crossing problem}: for a given survival
distribution
$p$ on $(0,\infty)$ find a function $b$ such that $b$
satisfies  (\ref{BCDefinition}), (\ref{pDef}). If such a boundary
exists and is unique, we
denote it by $\cB[p]$.
\end{enumerate}
The boundary crossing problem is classical, and the
subject of a large literature. The inverse boundary crossing problem
has recently
been the subject of increased interest by probabilists and researchers
in mathematical
finance. The main purpose of this paper is to show that the inverse boundary
crossing problem is well posed.

According to \citet{ZuccaSacerdote}, the problem
was originally posed by A. N. Shiryayev in 1976, for the special case where
$X_{t}$ is a Brownian motion and $p$ is the exponential distribution.
\citet{DudleyGutmann} and
\citet{Anulova} showed that there exists a \textit{stopping time} with
the given distribution; however,
this stopping time is not realized as the first time the process $X$
crosses a boundary $b$.
Recently, there has been an increase of interest in the problem due to
its importance in
applications. In mathematical finance, with $X_{t}$ an indicator of a
firm's financial
health, and $p$ the distribution of its time to default (estimated from
the prices of
market instruments), the problem is to find a default barrier that
reproduces the given
default distribution. Many authors have proposed numerical
methods for finding such a boundary, including \citet{HullWhiteTwo},
\citet{Is},
\citet{HuangTian}, \citet{AZ} and \citet{ZuccaSacerdote}. A formulation
of the
problem in terms of nonlinear Volterra integral equations has been given
by \citet{Pe1} [see also \citet{PeskirShiryayev} for a more detailed
discussion].

The numerical method proposed by \citet{AZ} is most relevant to our
work. They note that
for sufficiently smooth boundaries $b$, the function
$U(x,t) = \partial_x\mathbb{P}(\hat\tau\geq t, X_{t} \leq x)$
should be the solution of the free boundary
problem
%
%
\begin{equation}\label{UFBP}
\cases{
\cL_{1}U(x,t) = 0,
&\quad for $x>b(t)$, $t>0$,
\cr
U(x,t)=0, &\quad for $x\leq b(t)$, $t>0$,
\cr
U(x,0)=\rho_{0}(x,0), &\quad for $x\in\mathbb{R}$,}
\end{equation}
with the free boundary condition
%
%
\begin{equation}\label{pderivative}
\dot p(t)= -\tfrac1 2 (\sigma^2 U)_x
|_{x=b(t)}\qquad\forall t\geq0,
\end{equation}
where $\cL_{1}$ is the differential operator
%
%
\begin{equation} \label{1.cL1}
\cL_1\phi:=\frac{\partial\phi}{\partial t} -\frac12 \frac
{\partial^2[\sigma^2
\phi]}{\partial x^2}+ \frac{\partial[\mu\phi]}{\partial x}.
\end{equation}

\citet{AZ} perform a change
of variables to ``straighten out the boundary,'' and then solve the
resulting transformed PDE
numerically using finite differences.

An analytic study of the inverse boundary crossing problem was
initiated by \citet{CCCS}.
In that work, we defined
%
%
\begin{equation}
w(x,t) = \int_{x}^{\infty} U(y,t) \,dy.
\end{equation}
Formally, direct calculation from (\ref{UFBP}), (\ref{pderivative})
shows that $w$ should
satisfy the free boundary problem
%
%
\begin{equation} \label{WFBP}
\cases{\cL w(x,t) = 0, &\quad for $x>b(t)$, $t>0$,
\cr
w(x,t)=p(t), &\quad for $x\leq b(t)$, $t>0$,
\cr
w_{x}(x,t) = 0, &\quad for $x\leq b(t)$, $t>0$,
\cr
w(x,0)=1-p_{0}(x,0), &\quad for $x\in\mathbb{R}$,}
\end{equation}
where $\cL$ is the differential operator
%
%
\begin{equation}\label{1.cL}
\cL\phi: =\frac{\partial\phi}{\partial t} -\frac12
\frac{\partial}{\partial x}\biggl(\sigma^2\,\frac{\partial
\phi}{\partial x} \biggr)+ \mu\,\frac{\partial\phi}{\partial
x}.
\end{equation}
Formally,
%
%
\begin{equation}
w(x,t) = \mathbb{P}(\hat\tau\geq t, X_{t} > x),
\end{equation}
providing the connection between the probabilistic problems (\ref
{BCDefinition}),
(\ref{pDef}) and our analytic approach. Based on the free boundary
problem (\ref{WFBP}), one can infer that $w$ should satisfy the
variational inequality
%
%
\begin{eqnarray}
\label{2.var}
\max\{ \cL w , w-p \}&=& 0\qquad\mbox{in }L^\infty\bigl(\BbbR\times(0,\infty)\bigr),
\nonumber\\[-8pt]\\[-8pt]
w(\cdot,0)&=&1- p_0(\cdot,0) \qquad\mbox{on }\BbbR, \nonumber
\end{eqnarray}
and $b$ can be recovered from $w$ by
%
%
\begin{equation}
b(t) := \inf\{ x | p(t)-w(x,t)>0 \}\qquad\forall t>0.
\end{equation}
In \citet{CCCS}, the existence and
uniqueness of a viscosity solution to (\ref{2.var}) was proved.
However, no attempt
was made to connect the resulting functions $w,b$ to the original
probabilistic formulation
of the inverse boundary crossing problem. In this paper, we show that
$b$ does in fact
give a boundary that reproduces the survival distribution $p$. This is
complicated by the fact that it is
very difficult to prove the regularity of the boundary $b$ derived from
the variational
inequality.\footnote{The problem of the regularity of the boundary has
subsequently been
investigated by \citet{ChenRegularity}.}
As a consequence, in order to verify that $b$ has the required hitting
distribution,
we must first study the problem of computing the boundary crossing probabilities
of diffusions to nonsmooth boundaries. To this end, for a given function
$b\dvtx(0,\infty) \to[-\infty,\infty)$, we define
%
%
\begin{equation} \label{3a.bstar}\quad
b^*(t):=\max\Bigl\{b(t), \varlimsup_{s\to t}
b(s)\Bigr\},\qquad
b^*_-(t):=\varlimsup_{s\nearrow t} b(s) \qquad\forall t>0.
\end{equation}
When needed, we also define $b^*(0):=\varlimsup_{s\searrow0} b(s)$.
We also employ the notation $Q_b:=\{(x,t) | x>b(t),t>0\}$.
It will turn out that the inverse boundary crossing problem is most
naturally formulated
in the following spaces:
\begin{eqnarray*}
B_0&:=& \biggl\{ b\dvtx(0,\infty)\to[-\infty,\infty) |
b=b^*=b^*_-, \BbbP\biggl(\bigcup_{\vep>0} \bigcap_{s\in(0,\vep
)}\{
X_s\geq b(s)\}\biggr)=1 \biggr\},
\\
P_0 &:=& \{ p\in C([0,\infty)) |
p(0)=1\geq p(s)\geq p(t)>0\ \forall t> s\geq0\}.
\end{eqnarray*}
The main result of this paper is the following theorem.
\begin{theorem}\label{th3.1} \label{MainTheorem}

\begin{enumerate}
\item For every $p\in P_0$, there exists a unique viscosity solution,
$w$, for the survival distribution of the inverse boundary crossing problem
associated with $p$ [i.e., a~viscosity solution
of problem (\ref{2.var})]. In addition the unique solution, $w$,
satisfies
%
%
\begin{eqnarray}\label{monotonew}
0&\leq&1- [w(x,t)+p_0(x,t)]\leq1-p(t),\nonumber\\[-8pt]\\[-8pt]
w(x,t) &\leq& w(y,t)\qquad \forall t\geq0, x\in\BbbR, x\geq y.\nonumber
\end{eqnarray}
%

Consequently, the operator $\cB$
%
%
\begin{equation}
b(t) = \cB[p](t) :=\inf\{ x\in\BbbR| w(x,t)<p(t)\}\qquad
\forall t>0,
\end{equation}
is well defined on $P_0$.

\item For every $p\in P_0$, $\cB[p]\in B_0$ and $(\cP\circ
\cB)[p]=p$, where $(\cP[b])(t)=\mathbb{P}(\hat\tau\geq t)$ and
$\hat\tau$ is defined
as in (\ref{BCDefinition}).

\textup{(This implies that for a given $p\in P_0$, $b:=\cB[p]$ is a
solution of the inverse problem since
$\cP[b]=(\cP\circ\cB)[p]=p$.)}

\item For every $b\in B_0$, $\cP[b]\in P_0$ and $(\cB\circ
\cP)[b] =b$.

\textup{(This implies that for a given $p\in P_0$, if $\tilde b\in B_0$ is
a solution of the inverse problem, i.e., $\cP[\tilde b]=p$, then
$\tilde b$ can be identified as
$\tilde b=(\cB\circ\cP)[\tilde b]=\cB[p]$, the viscosity solution
of the inverse problem.)}

\item If $(w,b)$ is a classical \textup{(\textit{i.e.}, $w+p_0\in C(\BbbR\times
[0,\infty))$}, \textup{$\partial_x w\in C(\BbbR\times(0,\infty)$}, \textup{$\partial_t
w,\partial_{xx}w\in C(Q_b)$)} solution of the free
boundary problem (\ref{WFBP}), then $b$ is the solution of the
inverse boundary
crossing problem associated with $p$, that is, $\cP[b]=p$.
Similarly, if $(U,b)$ is a classical \textup{(\textit{i.e.}, $U-\rho_0\in C(\BbbR
\times[0,\infty))$}, \textup{$\partial_t U,\partial_{xx}^2 U \in C(Q_b)$)}
solution of the free
boundary problem (\ref{UFBP}), (\ref{pderivative}),
then $\cP[b]=p$.
\end{enumerate}
\end{theorem}

The proof of the above theorem proceeds as follows. We begin by
studying the direct
problem of computing the distribution of $\hat\tau$, and the function
$w(x,t) =\mathbb{P}(\hat\tau\geq t, X_{t} > x)$ for boundaries $b\in
B_{0}$. By considering
a carefully constructed discrete approximation scheme motivated by
(\ref{UFBP}) when $b$ is known, we are able to show that $w$
is the unique
viscosity solution to (\ref{WFBP}).
Elementary calculations verify that the viscosity solution of the
variational inequality (\ref{2.var}) also solves~(\ref{WFBP}). Once
we have also
determined that $\{x > b\} = \{ w < p \}$, the verification proceeds by
relatively straightforward
arguments.

We note that the sequence of stopping times constructed in our discrete
time approximation
actually converges to the first time that $X$ is strictly below the
boundary $b$, $\tau$ as
given by (\ref{TauDefinition}).
This definition of the boundary crossing time is slightly different
from the
standard one (\ref{BCDefinition}) for $\hat\tau$.
We have found that for the analytic approach we take here (particularly
for rough boundaries),
our definition is more convenient. In Section \ref{MeasureSection}
below, we show that for boundaries with minimal regularity properties
(including those arising from the solution to the inverse boundary
crossing problem, $b\in B_{0}$), $\mathbb{P}(\tau= \hat\tau)=1$.

%



The remainder of the paper is structured as follows. The second section
proves measurability
properties of $\tau$ and $\hat\tau$, and proves that these times are
almost surely equal.
In addition, it gives preliminary results that are needed for the study
of our approximation
scheme. The third section studies the approximation scheme in detail,
and proves
convergence.
The convergence provides a rigorous connection between the
probabilistic definition of the survival probability $p$ and the PDE
definition of the survival distribution $w$.
The fourth section formulates viscosity solutions for the \textit
{direct problem}
of computing $p$ for a given $b\in B_{0}$, and shows that the survival
distribution
$w$ gives the unique viscosity solution for the direct problem. The
fifth section
provides the link between the variational inequality studied in \citet
{CCCS} and the
inverse boundary crossing problem. It also provides a sufficient
condition under which
the resulting boundary $b$ is continuous.

\section{Crossing times of upper-semi-continuous boundaries}
\label{MeasureSection}
We calculate boundary crossing distributions for rough boundaries
based on discrete time approximations to be studied
in the next section. In order to ensure convergence of
our approximation scheme, the time points used must be chosen
carefully. We refer to the
selected points as the ``landmark points'' of the boundary. In this
section, we
begin by defining the landmark points and investigating their
properties. Using these properties,
we study the measurability of $\tau$ and $\hat\tau$, show that the
boundary crossing
times of $b$ and $b^{*}$ are equal and that $\mathbb{P}(\tau= \hat
\tau)=1$ for
$b\in B_{0}$.
\begin{definition}
Let $b\dvtx(0,\infty)\to[-\infty,\infty)$, and $b^{*}$ be its
upper-semi-continuous envelope. The set of
landmark points of $b$, denoted by $\LM(b)$, is defined as follows:
%
%
\begin{eqnarray}
\label{LMOne}
\LM(b)&:=&\bigcup_{n\in\BbbN}\LM_n(b),\qquad
\LM_n(b):=\{t_n^i | i\in\BbbN\},
\\
\label{3a.LM}
t_n^i &:=& \inf\biggl\{ t\in\biggl[ \frac{i}{2^{n}}, \frac
{i+1}{2^{n}}\biggr] \Big| b^*(t)\geq\sup_{[2^{-n}i,2^{-n}(i+1)]}
b(\cdot)\biggr\}.
\end{eqnarray}
\end{definition}

The following lemma summarizes some properties of the landmark points
that are used
in the paper.
\begin{lemma}\label{LMLemma} Let $b\dvtx(0,\infty)\to[-\infty,\infty
)$, and let its landmark
points $\LM(b)$ be defined as in (\ref{3a.LM}).
\begin{enumerate}
\item For $i,n\in\mathbb{N}$,
$b^{*}(t^{i}_{n}) \geq b^{*}(s)$
for every $s\in[2^{-n}i, 2^{-n}(i+1))$.
\item For $i,n\in\BbbN$, either $t_{n}^i=t_{n+1}^{2i}$ or
$t_n^i=t_{n+1}^{2i+1}$, so
$\LM_n(b)\subset\LM_{n+1}(b)$.
\end{enumerate}
\end{lemma}
\begin{pf}
1. If $s \in(2^{-n}i,2^{-n}(i+1))$, then by definition
$b^{*}(t^{i}_{n}) \geq b(t)$ for
$t\in(s-\varepsilon, s+\varepsilon)$ with $\varepsilon$ small enough.
$b^{*}(t^{i}_{n}) \geq b^{*}(s)$ follows immediately. If $s = 2^{-n}i$, and
$b^{*}(s) \geq b^{*}(t^{i}_{n}) \geq\sup\{b(t) | t\in[2^{-n}i,
2^{-n}(i+1)]\}$ then
$s =t^{i}_{n}$ by (\ref{3a.LM}).\vspace*{2pt}

2. We first claim that if $t^{i}_{n} \in
[2^{-n}i,2^{-(n+1)}(2i+1))$ then
$t^{i}_{n} = t^{2i}_{n+1}$. Clearly,
$b^{*}(t^{i}_{n}) \geq\sup\{b(s) | s\in[2^{-n}i,2^{-n}(i+1)]\}
\geq\sup\{ b(s)\vspace*{2pt} | s \in[2^{-n}i, 2^{-(n+1)}(2i+1)]\}$, so by definition
$t^{2i}_{n+1} \leq t^{i}_{n}$. If the inequality is strict, there is a
$\delta> 0$ small
enough so that $(t^{i}_{n}-\delta, t^{i}_{n}+\delta) \subseteq
(2^{-n}i, 2^{-(n+1)}(2i+1))$,
and since $b^{*}(t^{2i}_{n+1}) \geq b$ on this interval, we obtain
$b^{*}(t^{2i}_{n+1})
\geq b^{*}(t^{i}_{n})$, contradicting the definition of $t^{i}_{n}$. A
similar proof shows that
if $t^{i}_{n} \in[2^{-(n+1)}(2i+1),2^{-n}(i+1))$, then $t^{i}_{n} =
t^{2i+1}_{n+1}$. Finally,
it is easy to see that if $t^{i}_{n} = 2^{-n}(i+1)$, then
$\sup\{ b(s) | s\in[2^{-n}i,2^{-n}(i+1)]\}
= \sup\{ b(s) | s\in[2^{-(n+1)}(2i+1),2^{-n}(i+1)]\}
\geq\sup\{ b(s) | s\in[2^{-n}i,2^{-(n+1)}(2i+1)]\}$,
after which repeating the same argument by contradiction ensures that
$t^{i}_{n}=t^{2i+1}_{n+1}$.
\end{pf}

The following lemma collects some properties of upper-semi-continuous
functions that
are used throughout the paper. The proofs are elementary, and are omitted.
\begin{lemma}\label{USCLemma}
Let $b\dvtx(0,\infty)\to[-\infty,\infty)$ be upper-semi-continuous.
\begin{enumerate}
\item If $x\dvtx[0,\infty) \to(-\infty,\infty)$ is continuous, then for
all $t > 0$,
\[
\inf\{ s > 0 | x(s) \leq b(s)\} > t \quad\iff\quad x(s) > b(s) \qquad\forall s\in(0,t]
\]
\item The set $Q_b:=\{(x,t) | x>b(t),t>0\}$ is open.
\end{enumerate}
\end{lemma}

The next proposition addresses two main issues. First, it considers the
measurability
of $\tau$ and $\hat\tau$, to ensure that the various functions
considered in the remainder of the
paper are well defined. Second, it shows that for the purposes of
computing the distribution of
$\tau$, it is enough to consider the upper-semi-continuous envelope,
$b^{*}$, of the boundary $b$. We
observe that the result for $\tau$ holds with minimal assumptions on
the function $b$ (we have not
even assumed measurability). Furthermore, we note that some of the
results on measurability
could be derived by applying more general theorems [e.g., $\hat
\tau$ is the
first hitting time of the two-dimensional process $(X_{t},t)$ to the
set $\{ (x,s) | x \leq b(s) \}$,
which is closed when $b$ is upper-semi-continuous]. However,
we have decided to present elementary proofs of these assertions to
make the paper more
self-contained.
\begin{proposition}\label{newth1} Let $b\dvtx(0,\infty)\to[-\infty
,\infty)$.
\begin{enumerate}
\item Let $b^*$ be as in (\ref{3a.bstar}) and
$\LM(b)$ be as in (\ref{LMOne}) and (\ref{3a.LM}). Then
for every $t>0$,
%
%
\begin{equation}\label{3a.bt}
\bigcap_{s\in(0,t)} \{X_s \geq b(s)\}=
\bigcap_{s\in(0,t)\cap\LM(b)}
\{ X_s \geq b^*(s)\}.
\end{equation}

Consequently, we can define the \textup{first boundary crossing time}
$\tau\dvtx\Omega\to[0,\infty]$, the \textup{survival probability}
$p\dvtx[0,\infty)\to[0,1]$, and the
\textup{survival distribution} $w\dvtx\BbbR\times[0,\infty)\to[0,1]$ by
%
%
\begin{eqnarray}
\label{tauDefinition}
\tau(\omega) &:=& \inf\{ s>0 | X_s(\omega)< b(s)\}\qquad\forall
\omega\in\Omega,
\\
\label{pDefinition}
p(t) &:=& \BbbP\bigl(X_s \geq b(s)\ \forall s\in(0,t)\bigr) \qquad \forall t\geq0,
\\
\label{3a.pw}
w(x,t) &:=& \BbbP\bigl(X_s \geq b(s)\
\forall s\in(0,t), X_t >x\bigr)\qquad \forall x\in\BbbR, t\geq0.
\end{eqnarray}

In addition, $\tau$ is an optional time with respect to the filtration
generated by the process $X$,
$ \{\tau\geq t\} \in\cF^{X}_t, \forall t\geq0$. Also
\[
p(t)=\BbbP(\tau\geq t),\qquad
w(x,t)=\BbbP(\tau\geq t, X_t>x).
\]

\item Let the \textup{(conventional)} \textup{first crossing time} $\hat
\tau\dvtx\Omega\to[0,\infty]$ be defined by
\[
\hat\tau(\omega) :=\inf\{ s>0 | X_s(\omega)\leq b(s)\}\qquad
\forall\omega\in\Omega.
\]
If $b$ is upper-semi-continuous, i.e., $b=b^*$, then $\hat\tau$ is a
stopping time
with respect to the filtration generated by the process $X$,
$ \{\omega\in\Omega| \hat\tau(\omega) >t\} \in\cF^{X}_t\
\forall t\geq0$, so that
we can define
$\hat p(t):=\BbbP(\hat\tau(\omega)>t)\ \forall t\geq0$.
\end{enumerate}
\end{proposition}
%
%
\begin{pf}
1.
When $t=0$, we have $\{\tau\geq0\} =\Omega=\{\omega\in\Omega|
X_s(\omega)
\geq b(s), \operatorname{Con}$ $\forall s\in(0,t)\}$ and $p(0)=\BbbP(\Omega)=1$.
Now we assume that $t>0$.
It is easy to verify that
\[
\bigcap_{s\in(0,t)} \{ X_s\geq b(s)\}=\{X_s\geq
b(s)\ \forall
s\in(0,t)\}
= \{\tau\geq t\}.
\]
Hence, to complete the proof of the first assertion, it suffices to
verify (\ref{3a.bt}).
Suppose $\omega\in\{X_s \geq b(s)\ \forall s\in(0,t)\}$. Then
$X_s(\omega)\geq b(s)$ for every $s\in(0,t)$. For every $\hat s\in
(0,t)$, by the continuity of
$X_\bullet(\omega)$,
\[
X_{\hat s}(\omega) = \lim_{s\to\hat s} X_{s}(\omega
) \geq
\max\Bigl\{ b(\hat s), \varlimsup_{s\to\hat s} b(s)\Bigr\}
=b^*(\hat s).
\]
As $\hat s\in(0,t)$ is arbitrary, we have $\omega\in\{X_s\geq b^*\
\forall s\in(0,t)\}$.
Thus, $ \{X_s\geq b(s)\ \forall s\in(0,t)\} \subset\{X_s\geq
b^*(s)\ \forall s\in(0,t)\}\subset\{X_s\geq b^*(s)\ \forall
s\in(0,t)\cap\LM(b)\}$.

Next, suppose $\omega\in\{X_s(\omega)\geq b^*(s) \ \forall s\in
(0,t)\cap\LM(b)\}$. Let $\hat s\in(0,t)$ be arbitrary. We want to
show that $X_{\hat s}\geq b(\hat s)$. For each integer\vspace*{1pt} $n$ satisfying
$2^{-n}\leq\hat s$, let
$i_n$ be the integer such that $\hat s\in[i_n 2^{-n},(i_n+1)2^{-n})$.
Then $t_n^{i_n}\in\LM(b)$ and $X_{t_n^{i_n}}(\omega)\geq
b^*(t_n^{i_n})\geq b^*(\hat s)$
by Lemma \ref{LMLemma}. Hence,
\[
X_{\hat s}(\omega) = \lim_{n\to\infty} X_{t_n^{i_n}}(\omega) \geq
\varlimsup_{n\to\infty} b^*(t_n^{i_n}) \geq b^*(\hat s)\geq b(\hat s).
\]
Since $\hat s$ is arbitrary, we see that
$\omega\in\{X_s\geq b(s)\ \forall s\in(0,t)\}$. Consequently,
\begin{eqnarray*}
\{\tau\geq t\}&=&\bigcap_{s\in(0,t)} \{X_s\geq b(s)\}
=\bigcap
_{s\in(0,t)}\{X_s\geq b^*(s)\}\\
&=&\bigcap_{s\in(0,t)\cap\LM(b)} \{
X_s\geq b^*(s)\}\in\cF^{X}_t.
\end{eqnarray*}
This proves (\ref{3a.bt}) and also the first assertion.

2. Assume that $b$ is usc, that is, $b=b^*$. Then by the
continuity of the sample paths of $X$,
if $\hat s>0$ and $X_{\hat s}(\omega)>b(\hat s)$ then there exists
$\delta>0$ such that
$X_s(\omega)>b(s)+\delta$ for all $s\in[\hat s-\delta,\hat s+\delta
]$. By the Heine--Borel theorem,
if $X_s>b(s)$ for every $s\in[a,c]\subset(0,\infty)$, then there exists
a large integer $i$ such that
$X_s>b(s)+2^{-i}$ for every $s\in[a,c]$. Hence,
for every $t>0$,
%
%
\begin{eqnarray} \label{3b.bt}
\{\hat\tau>t\}&=&\bigcap_{s\in(0,t]} \{X_s>b(s)\} =
\bigcap
_{n\in\BbbN} \bigcap_{s\in[2^{-n}t,t]}\{X_s>b(s)\}\nonumber\\
&=&\bigcap_{n\in\BbbN} \bigcup_{i\in\BbbN} \bigcap
_{s\in[2^{-n}t,t]} \{X_s\geq b(s)+2^{-i}\}
\\ &=& \bigcap_{n\in\BbbN} \bigcup_{i\in\BbbN}
\bigcap_{s\in\{2^{-n}t,t\}\cup(\LM(b)\cap[2^{-n}t,t])}
\{X_s\geq b(s)+2^{-i}\}\in\cF_t^{X}.\nonumber
\end{eqnarray}
This completes the proof of the second assertion.
\end{pf}

The following proposition justifies our choice to work with $\tau$,
the first time the process is
strictly below the boundary, rather than $\hat\tau$, the first time
the process hits the boundary.
In particular, the second assertion implies that these times are almost
surely equal, and hence
they have the same distributions (so solving the inverse boundary
crossing problem for $\tau$
is the same as solving it for $\hat\tau$). We will see in the next
section that it is easier to work with
$\tau$ in defining approximations for the boundary crossing problem.

\begin{proposition} \label{ASProposition}
Let $b\dvtx(0,\infty)\to[-\infty,\infty)$.
\begin{enumerate}
\item Let $b^*$ and $b^*_-$ be as in (\ref{3a.bstar}). Then
\[
\{\omega\in\Omega| \tau(\omega)=t \mbox{ or }\hat\tau
(\omega)=t\} \subset\{\omega\in\Omega| X_t(\omega)\in
[b^*_-(t),b^*(t)]\}\qquad\forall t>0.
\]

\item The set $\{\tau\not=\hat\tau\}$
has
probability zero, so that
\begin{eqnarray*}
p(t)=p(t-)=\hat p(t-) &=& \BbbP(\tau\geq t)
\qquad\forall t>0,
\\
\hat p(t)=\hat p(t+)=p(t+) &=& \BbbP(\tau>t) \qquad\forall t\geq0.
\end{eqnarray*}
Consequently,
if $b^*=b^*_-$, then $p\in C((0,\infty))$, $\hat p=p$ on $(0,\infty
)$, and $\hat p\in C([0,\infty))$.
\end{enumerate}
\end{proposition}
\begin{pf}
To prove the first assertion, let $t>0$ and $\omega\in\{\tau=t\}\cup
\{\hat\tau=t\}$.
Then $X_s(\omega)\geq b(s)$ for all $s\in(0,t)$ so
$ X_t(\omega) \geq\varlimsup_{s\nearrow t} b(s)=b_-^*(t). $
Also, $\tau(\omega)=t$ or $\hat\tau(\omega)=t$ implies that there
exists a sequence $\{s_i\}$ of positive numbers such that $\lim_{i\to
\infty} s_i=t$ and $X_{s_i}(\omega)\leq b(s_i)$ for all $i$. This
implies that
$ X_t(\omega) =\lim_{i\to\infty} X_{s_i}(\omega)\leq\varlimsup
_{i\to\infty} b(s_i) \leq b^*(t). $
Hence, $X_t(\omega)\in[b_-^*(t),b^*(t)]$.
Also, note that if $b_-^*(t)=b^*(t)$, then $\BbbP(\{\tau=t\}\cup\{
\hat\tau=t\})\leq
\BbbP(\{X_t=b^*(t)\})\,{=}\,0$.

Since the family $\{\tau>t+\vep\}_{\vep\geq-t}$ of sets is
monotonic in $\vep$, we see that
\begin{eqnarray*}
p(t+)&=&\lim_{\vep\searrow0} p(t+\vep) =\lim_{\vep\searrow0}\BbbP
(\tau\geq t+\vep)
\\
&=&\BbbP\biggl(\bigcup_{\vep>0}\{\tau\geq t+\vep\}\biggr)=\BbbP(\tau>t)\qquad
\forall t\geq0,
\\
p(t-)&=&\lim_{\vep\searrow0} p(t-\vep)=\lim_{\vep\searrow0}\BbbP
(\tau\geq t-\vep)
\\
&=&\BbbP\biggl(\bigcap_{\vep>0}\{\tau\geq t-\vep\}\biggr)
=\BbbP(\tau\geq t)=p(t)\qquad\forall t>0.
\end{eqnarray*}

Clearly, if $b_-^*(t)=b^*(t)$, then $p(t-)-p(t+)=\BbbP(\tau=t)=0$ so
$p$ is continuous at $t$.
Similarly, when $b^*=b$ so $\hat p$ is well defined, we have $\hat
p(t)=\BbbP(\hat\tau>t)=\hat p(t+)$ and $\hat p(t-)=\BbbP(\hat\tau
\geq t)$.

To complete the proof,
it remains to show that the set $\{\tau\not=\hat\tau\}$ has
probability zero.
For every $\omega\in\Omega$,
\[
\tau(\omega)=\inf\{s>0 | X_s(\omega)<b(s)\} \geq \inf\{
s>0 | X_s(\omega)\leq b(s)\}=\hat\tau(\omega).
\]

Now, suppose $\tau(\omega)\not=\hat\tau(\omega)$. Then we must
have $ \hat\tau(\omega)<\tau(\omega)$. Set $t=\tau(\omega)$.
Then $X_s\geq b(s)\ \forall s\in(0,t)$. By continuity, we also have
$X_s\geq b^*(s) \ \forall s\in(0,t)$.
Set $\hat t=\hat\tau(\omega)$.
If $\hat t=0$, then by definition, there exists $r\in(0,t)$ such that
$X_r(\omega)\leq b(r)$.
If $\hat t>0$, then by definition, there exists $r\in[\hat t,t)$ such
that $X_r(\omega)\leq b(r)$.
As $X_s(\omega)\geq b^*(s)$ for all $s\in(0,t)$, in either case, we
have $r\in(0,t)$ and $X_r(\omega)=b^*(r)$. Taking $r_1\in\LM(b)\cap
(0,r)$ and $r_2\in\LM(b)\cap(r,t)$ we obtain
$\min_{s\in[r_1,r_2]} \{X_{s}(\omega)-b^*(s)\}=0$. Hence,
\[
\{\tau\not=\hat\tau\} =\{\hat\tau<\tau\} \subset\bigcup
_{r_1\in\LM(b)}
\bigcup_{r_2\in\LM(b)\cap(r_1,\infty)} \bB(r_1,r_2),
\]
where for every $0<a<c<\infty$,
\[
\bB(a,c)
= \Bigl\{\omega\in\Omega\big| \min_{s\in[a,c]} \{X_s(\omega)-b^*(s)\}
=0\Bigr\}.
\]

Note that for each $c>a>0$, ${\mathbf B}(a,c)$ is $\cF_c$ measurable since
\begin{eqnarray*}
{\mathbf B(a,c)} &=& \Bigl\{\min_{s\in[a,c]} \{
X_s-b^*_s\}\geq
0\Bigr\}\\
&&{}\Bigm\backslash
\bigcup_{n=1}^\infty\Bigl\{\min_{s\in[a,c]} \{X_s -b^*(s)\}\geq
2^{-n}\Bigr\},
\\
\Bigl\{\min_{s\in[a,c]} \{X_s-b^*(s)\}\geq h\Bigr\}
&=& \bigcap_{s\in\{a,c\}\cup(\LM(b)\cap[a,c])} \{X_s\geq b^*(s)+h\}
\in\cF_c \qquad\forall h\in\BbbR.
\end{eqnarray*}

It remains to show that for each $c>a>0$, the set
$\bB(a,c)$ has measure zero. Suppose, on the contrary, that
$\BbbP(\bB(a,c))>0$ for some fixed $c>a>0$. Fixing $t_0\in(0,a)$, we
then have
\[
0<\BbbP({\mathbf B}(a,c) ) =\int_\BbbR\BbbP\bigl( \bB
(a,c)|X_{t_0}=z\bigr)\rho_0(z,t_0) \,dz.
\]
Consequently, there exists a finite number $M>0$ such that
\[
\int_{-M}^M \BbbP\bigl( \bB(a,c)|X_{t_0}=z\bigr)\rho_0(z,t_0)
\,dz>0.
\]

For each $h\in\BbbR$, we consider the set
\[
\bB^h(a,c) = \Bigl\{ \min_{s\in[a,c]} \bigl(X_s-b^*(s)\bigr)=h\Bigr\}.
\]
For the process $\{X_t\}_{t\geq t_0}$,
for each $\omega\in\{X_s=z\}$ and $h\in\BbbR$, we denote by
$\omega^h$ the element in $\{X_{t_0}=z+h\}$ such that
$ X_t(\omega^h)= h+ X_
t(\omega)\ \forall t\in[t_0,\infty)$.
Then
\[
\min_{s\in[a,c]} \bigl(X_s(\omega)-b^*(s)\bigr)=0
\quad\Longleftrightarrow\quad \min_{s\in[a,c]} \bigl(X_s(\omega^h)-b^*(s)\bigr)=h.
\]
Assume for simplicity that we are dealing with Brownian motion. [By a
change of variables,
we can assume $\sigma\equiv1$; see Section \ref{sec4}. Since $\mu$ is smooth
and bounded,
if $\{X_t\}$ is not a Brownian motion, we can use the Girsanov
theorem [\citet{KaratzasShreve}] to change to an equivalent
measure under which $X_{t}$ is a Brownian motion, and the argument
below can still be
used to show that $\mathbb{P}(\bB(a,c))=0$.]
By the translation invariance of Brownian motion
and the Markov property, we have
\[
\BbbP\bigl( \bB(a,c)|X_{t_0}=z\bigr)= \BbbP\bigl(\bB
^h(a,c)|X_{t_0}=z+h\bigr).
\]

Hence,
\begin{eqnarray*}
\BbbP(\bB^h(a,c)) &=& \int_\BbbR\BbbP\bigl( \bB
^h(a,c)|X_{t_0}=z+h\bigr)\rho_0(z+h,t_0) \,dz
\\
&=& \int_\BbbR\BbbP\bigl( \bB(a,c)|X_{t_0}=z\bigr)\rho_0(z+h,t_0) \,dz
\\
&\geq& \min_{z\in[-M,M]}\frac{\rho_0(z+h,t_0)}{\rho_0(z,t_0)}
\int_{-M}^M \BbbP\bigl( \bB(a,c)|X_{t_0}=z\bigr)\rho_0(z,t_0) \,dz>0.
\end{eqnarray*}

Note that all elements in $\{\bB^h(a,c)\}_{h\in\BbbR}$ are disjoint
and measurable.
We then obtain a contradiction since $\Omega$ does not contain an
uncountable disjoint union of measurable sets with positive
probability. Thus, $\bB(a,c)$ must have
probability zero, for every pair $(a,c)$ with $a>c>0$. Consequently,
the set $\{\hat\tau\not=\tau\}$ has probability zero. This
completes the proof of Proposition \ref{ASProposition}.
\end{pf}

\begin{theorem}
The operator $\cP$ defined by
$\cP[b](t) = \BbbP(\tau\geq t)$
maps $B_0$ to~$P_0$.
\end{theorem}

\begin{pf}
Suppose $b\in B_0$.
Then $b^*=b^*_-$ so $p:=\cP[b]\in C((0,\infty))$. In addition,
\begin{eqnarray*}
\lim_{\vep\searrow0} p(\vep) &=& \lim_{\vep\searrow0} \BbbP\biggl(\bigcap
_{s\in(0,\vep)} \{ X_s>b(s)\}\biggr)
\\
&=& \BbbP\biggl(\bigcup_{\vep>0} \bigcap_{s\in(0,\vep)}
\{X_s>b(s)\}\biggr) =1=p(0).
\end{eqnarray*}

Hence, $p\in C([0,\infty))$.
It remains to show that $p>0$ on $[0,\infty)$.
Since $p(0+)=p(0)=1$, there exists $\vep>0$ such that $p(t)>0$ for
every $t\in[0,\vep]$.
Let $T>\vep$. The upper-semi-continuity of $b$ implies that
$M:=\sup_{s\in[0,T]} b(s)$ is finite.
Then $\BbbP(\tau>\vep, X_\vep>M)>0$.
Using standard results for a
constant barrier $\tilde b\equiv M$ on the set $\{\tau>\vep,X_\vep
>M\}$ for the time interval $[\vep,T]$, we see that
$\BbbP(\{\tau>\vep, X_s>M$ $\forall s\in[\vep,T]\})>0$. Hence,
$p(T)>0$. As $T$ is arbitrary,
we see that $p>0$ on $[0,\infty)$, so that $p\in P_0$.
\end{pf}

\begin{proposition}[(A semi-continuous dependence
property of $\cP$)] \label{prop4}
Assume that $b,b_1,b_2,\ldots$ are
upper-semi-continuous functions having the property
\[
b_1\leq b_{2}\leq b_3\leq\cdots,\qquad b(t)=\lim
_{n\to
\infty} b_n(t)\qquad \forall t>0.
\]
Let $p=\cP[b]$ and $p_n=\cP[b_n]$. Then for every $t\geq0$,
$p(t)=\lim_{n\to\infty} p_n(t)$.
\end{proposition}
\begin{pf}
For every $t>0$,
\begin{eqnarray*}
\{\tau\geq t \} &=& \bigcap_{s\in(0,t)}\{ X_s\geq
b(s)\}
= \bigcap_{s\in(0,t)} \bigcap_{n\in\BbbN} \{X_s\geq b_n(s)\}\\
&=&
\bigcap_{n\in\BbbN} \bigcap_{s\in(0,t)} \{ X_s\geq b_n(s)\}.
\end{eqnarray*}
Hence,
\begin{eqnarray*}
p(t) &=& \BbbP\biggl(
\bigcap_{n\in\BbbN} \bigcap_{s\in(0,t)} \{ X_s\geq b_n(s)\}\biggr)
=\lim_{n\to\infty}\BbbP\biggl(\bigcap_{s\in(0,t)} \{ X_s\geq b_n(s)
\}\biggr) \\
&=&\lim_{n\to\infty} p_n(t).
\end{eqnarray*}
\upqed\end{pf}

\section{Approximating sequences for boundary crossing times}
In this section, we use the landmark points to construct
straightforward approximations that eventually, upon passing to the
limit, will allow us to transfer the problem of calculating the
survival probability to problems of solving partial differential equations.
The advantage of studying the first time the process is strictly below
the boundary
is suggested by comparing the relative complexity of the expressions
(\ref{3a.bt})
and (\ref{3b.bt}). In this case, a simple
approximation to the survival probability and distribution can be developed.
We approximate a real barrier $b$ by a simple barrier $b_n$
defined by $b_n(t)=b^*(t)$ if $t\in\LM_n(b)$ and $b_n(t)=-\infty$
otherwise. The approximate problem then involves only the random
variables $\{X_t | t\in\LM_n(b)\}$ so that all relevant
probabilities can be calculated through transition probability
densities. Though it turns out that survival probabilities computed
using both viewpoints are equivalent, we do not see a simple adaptation
of our
method that allows us to approximate $\hat\tau$ directly without
appealing to the results in
Section \ref{MeasureSection}.
\begin{proposition} \label{prop5} Let $b\dvtx[0,\infty)\to[-\infty
,\infty)$ be
upper-semi-continuous and $(\tau,p,w)$ be defined as in
(\ref{tauDefinition})--(\ref{3a.pw}). Let $\LM(b)=\bigcup
_{n\in\mathbb{N}} \LM_{n}(b)$
be the landmark points of $b$, $\LM_{n}(b)=\{t_n^i | i\in\BbbN\}$ and
%
%
\begin{eqnarray}
\tau_n(\omega) &:=& \min\{ s\in\LM_n(b) | X_s(\omega)< b(s)\},\qquad
p_n(t) := \BbbP(\tau_n\geq t),
\\
\label{2.wn}\quad
w_n(x,t) &:=& \BbbP(\tau_n\geq t, X_t>x).
\end{eqnarray}


Then the following hold:

\begin{enumerate}
\item For all $(x,t)$,
\[
w_{n}(x,t) = \int_{x}^{\infty} U_{n}(y,t) \,dy,
\]
where, for $t\in[0,t_n^0], \{\tau_n\geq t\}=\Omega, p_n(t)=1,
w_n(\cdot,t)=1-p_0(\cdot,t)$ and $U_n(\cdot,t)=\rho_0(\cdot,t)$, and
when $t\in(t_n^k,t_n^{k+1}]$ with $k\in\BbbN$, $\{\tau_n\geq t\}
=\bigcap_{i=1}^k \{X_{t_n^i}\geq b(t_n^i)\}$ and
%
%
\begin{eqnarray}\label{2.Un}
U_n(x,t)&=& \int_{b(t_n^k)}^\infty U_n(y,t_n^k) \rho(y,t_n^k;
x,t) \,dy,
\nonumber\\[-8pt]\\[-8pt]
p_n(t)&=&\int_\BbbR U_n(y,t)\,dy =\int
_{b(t_n^k)}^\infty U(y,t_n^k) \,dy.\nonumber
\end{eqnarray}



\item For every $n\in\mathbb{N}$,
$ \tau_n \geq\tau_{n+1}\geq\tau,
p_n \geq p_{n+1}\geq p, w_n\geq w_{n+1}\geq w,
\rho_0\geq U_n \geq U_{n+1} \geq0$.

\item There exists
$U\dvtx\BbbR\times(0,\infty)\to[0,\infty)$ such that for every $
\omega\in\Omega,t>0$ and $x\in\BbbR$,
\[
\lim_{n\to\infty}
(\tau_n(\omega),
p_n(t), w_n(x,t) ,U_n(x,t))=(\tau(\omega
),p(t),w(x,t),U(x,t)).
\]
\end{enumerate}
\end{proposition}
\begin{pf}
1. This result follows immediately from the Chapman--Kolmogorov
equations and
the fact that the fundamental solution of $\mathcal{L}_{1}$ gives the
transition
densities of the Markov process $X$ [see, e.g.,
\citet{FriedmanSDE}, Theorem I.6.5.4, page 149].
From the definition of $\tau_n$ and $\LM_n(b)$, it is easy to see that
$\{\tau_n\geq t_n^1\}=\Omega$ and $\{\tau_n\geq t\}=\bigcap_{i=1}^k\{
X_{t_n^k}\geq b(t_n^i)\}$ when $t\in(t_n^k,t_n^{k+1}]$.
When $t\in[0,t_n^1]$, $\{\tau_n\geq t\}=\Omega$ so the evaluation of
$p_n,w_n,U_n$ is trivial.
When $t\in(t_n^i,t_n^{i+1}]$, $\BbbP(\tau_n\geq t,X_t>x)=
\BbbP(X_{t_n^1}\geq b(t_n^1),\ldots,X_{t_n^i}\geq b(t_n^i),X_t>x)$,
so using the transition probability density for the Markov process, we have
$ U_n(\cdot,t) = \int_{b_n(t_n^i)}^\infty U_n(y,t_n^i) \rho(y,t_n^i;
\cdot,t)\,dy$,
from which we find the corresponding $w_n$ and $p_n$. The first
assertion thus follows.

2. By the second part of Lemma \ref{LMLemma},
we have $\tau\leq\tau_{n+1}\leq\tau_n$, and therefore
$p\leq p_{n+1}\leq p_n$, and $w\leq w_{n+1}\leq w_n$.
It is clear from the definition that $t_{n}^{0} \leq t_{n+1}^{0}$ and so
$\rho_{0} = U_{n} = U_{n+1}$ on $(0,t_{n+1}^{0}]$. Now suppose
$U_{n+1} \leq U_{n}$
on $(0,t_{n+1}^{k}]$. Let $t\in(t_{n+1}^{k},t_{n+1}^{k+1}]$. Then
$t\in(t_{n}^{j},t_{n}^{j+1}]$
for some $j$ (the case $t \leq t_{n}^{1}$ is easier and handled
similarly). Then
\begin{eqnarray*}
U_{n+1}(x,t) &=& \int_{b(t_{n+1}^{k})}^{\infty} U_{n+1}(y,t_{n+1}^{k})
\rho(y,t_{n+1}^{k};x,t) \,dy \\
&\leq&\int_{b(t_{n+1}^{k})}^{\infty} U_{n}(y,t_{n+1}^{k})
\rho(y,t_{n}^{k};x,t) \,dy \\
&=& \int_{b(t_{n+1}^{k})}^{\infty} \rho(y,t_{n+1}^{k};x,t) \int
_{b(t_{n}^{j})}^{\infty}
U_{n}(z,t_{n}^{j}) \rho(z,t_{n}^{j} ;y,t_{n+1}^{k}) \,dz \,dy
\\
&\leq&
\int_{b(t_{n}^{j})}^{\infty} U_{n}(z,t_{n}^{j}) \int_{-\infty
}^{\infty}
\rho(z,t_{n}^{j};y,t_{n+1}^{k})\rho(y,t_{n+1}^{k};x,t) \,dy \,dz
\\
&=& \int_{b(t_{n}^{j})}^{\infty} U_{n}(z,t_{n}^{j}) \rho
(z,t_{n}^{j};x,t) \,dz = U_{n}(x,t)
\end{eqnarray*}
%
and $U_{n+1}\leq U_{n}$ by induction on $k$.
The proof that $0\leq U_{n} \leq\rho_{0}$ is similar. The second
assertion thus follows.

3. The monotonicity of $(\tau_n,p_n,w_n,U_n)$ implies the
existence of the limit as $n\to\infty$.
First we show that $\lim_{n\to\infty}\tau_n=\tau$. For this, let
$\omega\in\Omega$ be arbitrary.
(i) If $\tau(\omega)=\infty$, then we have $\tau_n(\omega)=\infty
$ for all $n\in\BbbN$ so
$\lim_{n\to\infty}\tau_n(\omega)=\infty=\tau(\omega)$. (ii)
Suppose $\tau(\omega)<\infty$.
Note that (\ref{3a.bt}) gives
%
%
\begin{eqnarray}
\{\tau\geq t\} &=& \bigcap_{s\in\LM(b)\cap(0,t)}\{ X_s\geq b(s)\}
\nonumber\\[-8pt]\\[-8pt]
&=&\bigcap_{n\in\BbbN} \bigcap_{s\in\LM_n(b)\cap(0,t)} \{X_s\geq
b(s)\}
=\bigcap_{n\in\BbbN} \{ \tau_n\geq t\}.\nonumber
\end{eqnarray}
Set $t:=\lim_{n\to\infty}\tau_n(\omega)$. Then as $\tau_{n+1}\geq
\tau_n\geq\tau$ for all $n\in\BbbN$, we see that $\tau_n(\omega
)\geq t\geq\tau(\omega)$ for every $n\in\BbbN$. Consequently,
$\omega\in\bigcap_{n\in\BbbN}\{\tau_n\geq t\}=\{\tau\geq t\}$.
Hence, we must have $\tau(\omega)=t=\lim_{n\to\infty}\tau
_n(\omega)$. Combining the two cases, we obtain $\lim_{n\to\infty}
\tau_n(\omega)= \tau(\omega)$ for every $\omega\in\Omega$.

Next, we consider the limits of $w_n$ and $p_n$. When $t=0$, we have
$w(\cdot,0)=w_n(\cdot,0)=1-p_0(\cdot,0)$ and $p(0)=1=p_n(0)$.
When $t>0$, for each $x\in\BbbR$,
\begin{eqnarray*}
w_n(x,t)-w(x,t)&=&\BbbP(\tau_n\geq t,X_t>x)-\BbbP
(\tau\geq t,X_t>x)
\\ &=& \BbbP(\tau_n\geq t>\tau,X_t>x)\leq\BbbP(\tau<t\leq\tau
_n).
\end{eqnarray*}
Thus,
%
%
\begin{eqnarray}
&&
{\varlimsup_{n\to\infty}}\|w_n(\cdot,t)-w(\cdot,t)\|_{L^\infty
(\BbbR)}\nonumber\\
&&\qquad\leq{\varlimsup_{n\to\infty}} |p_n(t)-p(t)|
= \varlimsup_{n\to\infty} \BbbP(\tau_n\geq\tau>t)
\\
&&\qquad=\BbbP\biggl(\bigcap_{n\in\BbbN}\{\tau_n\geq t>\tau\}\biggr)=0.
\nonumber
\end{eqnarray}


Finally, defining $U:=\lim_{n\to\infty}U_n$ we complete the proof of
the proposition.
\end{pf}

The approximating functions $U_n$ introduced in the previous proposition
are expressed in terms of the transition densities of the diffusion
$X$. From
an analytic point of view, they are obtained step by step, for
$i=1,2,\ldots,$
by \mbox{solving} the diffusion equations $\cL_1 U_n=0$ in the set $\BbbR
\times(t_{n}^i,t_n^{i+1}]$ with initial values
$U(\cdot,t_n^i+)=U(\cdot,t)\cdot\chi_{(b(t_n^i),\infty)}(\cdot)$, where
$\chi_A(x)$ is the indicator function of the set $A$.
In the sequel, $\cL$ and $\cL_1$ are the differential operators
introduced in (\ref{1.cL}) and
(\ref{1.cL1}), respectively.
Recall the notation $Q_b:=\{(x,t) | x>b(t),t>0\}$.
When $b\dvtx[0,\infty)\to[-\infty,\infty)$ is upper-semi-continuous,
the set $Q_b$ is an open set with $[b(0),\infty)\times\{0\}$ as its
``initial'' boundary.
%
%
\begin{proposition}\label{newth2} Let $b\dvtx[0,\infty)\to[-\infty
,\infty)$ be upper-semi-continuous
and
$(p,w)$ be the survival probability and survival distribution
associated with $b$, defined in
(\ref{tauDefinition})--(\ref{3a.pw}).
Then there exists a function $U$ such that the following hold:
%
%
\begin{eqnarray}
\label{3a.pwU}
p(t)&=&\int_{-\infty}^\infty U(y,t) \,dy = \int
_{b(t)}^\infty U(y,t) \,dy, \\
w(x,t)&=&\int_x^\infty U(y,t) \,dy
\\[2pt]
\label{3a.ini}\\\nonumber\\[-35pt]
&&\eqntext{\forall x\in\BbbR, t>0,
0\leq U\leq\rho_0, 0\leq w\leq1-p_0,}\nonumber\\[-25pt]
\end{eqnarray}
\begin{eqnarray}\nonumber\\[-30pt]
&\|w(\cdot,t)+p_0(\cdot,t)-1\|_{L^\infty(\BbbR)} \leq1-p(t)\qquad
\forall t>0,&\\
\label{3a.eqUw}
&\cL_1 U \leq 0, \qquad\cL w\leq0\qquad
\mbox{in }\BbbR\times(0,\infty),&\nonumber\\[-8pt]\\[-8pt]
&\hspace*{-35pt}\cL_1 U=0,\qquad \cL w=0 \qquad\mbox{in }Q_b,&\nonumber
\end{eqnarray}
where the inequalities in (\ref{3a.eqUw}) are understood in the sense
of distributions.
\end{proposition}
\begin{pf}
Let $U=\lim_{n\to\infty} U_n$.
Since $\rho_0\geq U_n\geq U_{n+1}\geq0$, using the Dominated Convergence
theorem and the identity $w_n(x,t)=\int_x^\infty U_n(y,t) \,dy$
we obtain $w(x,t)=\int_x^\infty U(y,t)\,dy$ for every $x\in\BbbR$ and $t>0$.
Since $\tau(\omega)\geq t$ implies $X_t(\omega)\geq b(t)$, we see that
$w(x,t)=w(-\infty,t)=p(t)$ for every $x<b(t)$. 

It is clear that $0\leq U\leq\rho_0$ and $0\leq w\leq1-p_0$. Also,
for $t>0$,
\begin{eqnarray*}
w(x,t)&=&\BbbP(X_t>x)-\BbbP(\tau<t,X_t>x)
\geq\BbbP(X_t>x)-\BbbP(\tau<t)\\
&=&[1-p_0(x,t)]-[1-p(t)].
\end{eqnarray*}
Thus, $0\leq1-w(x,t)-p_0(x,t) \leq1-p(t)$ or $\|w(\cdot,t)+p_0(\cdot
,t)-1\|_{L^\infty(\BbbR)} \leq1-p(t)$. 


It is useful to note that for each $t>0$, both $w_n(\cdot,t)$ and
$U_n(\cdot,t)$ are smooth functions. In
addition, as functions of $(x,t)$, $w_n$ and $U_n$ are smooth in
$\BbbR\times(0,\infty)\setminus\bigcup_{i=1}^\infty
(-\infty,b_{n}(t_n^i)]\times\{t_n^i\}$.
In particular,
\[
\cL_1 U_n = 0,\qquad \cL w_n=0 \qquad\mbox{in }
Q_{b}:=\{(x,t) | x>b(t), t>0\}.
\]
Since both $\{U_n\}_{n\in\BbbN}$ and $\{w_n\}_{n\in\BbbN}$ are
uniformly bounded in any compact
subset of~$Q_b$, it then follows from standard results on parabolic
partial differential equations [see \citet{FriedmanParabolic}, Theorem
3.11, page 74, and
Theorem 3.15, page 80]
that $w,U\in C^\infty(Q_b)$ and
$\cL w=0$ and $\cL_1 U=0$ in $Q_b$.

The set of discontinuities of $U_n$ and $w_n$ is $\bigcup_{i\in\BbbN}
(-\infty,b(t_n^i)]\times\{t_n^i\}$. In particular,
\begin{eqnarray*}
w_n(\cdot, t)&=&w_n(\cdot,t-),\qquad
U_n(\cdot,t)=U_n(\cdot,t-),\\
p(t)&=&p(t-) \qquad\forall t \in(0,\infty),
\\
U_n(x,t+)&=&0,\qquad
w_n(x,t+)=p_n(t) \qquad\forall x<b(t), t\in
\{t_n^i\}_{i\in\BbbN}.
\end{eqnarray*}
Denote by $\delta(\cdot-s)$ the Dirac measure with mass at $s$ and by
$\chi_{A}$ the characteristic function of the set $A$. Then in the
sense of distributions, we find that
\begin{eqnarray*}
\cL_1 U_n &=& \sum_{i=0}^{\infty} [ U_{n}(x,\cdot)
]_{t_{n}^{i} -}^{t_{n}^{i} +}
\delta(t-t_{n}^{i})
\\
&=& -\sum_{i=0}^\infty U_n(x,t_n^i) \delta(t-t_n^i) \chi
_{(-\infty,b(t_n^i)]}(x)\leq0\qquad
\mbox{in }\BbbR\times(0,\infty).
\\
\cL w_n &=& \sum_{i=0}^{\infty}[w_{n}(x,\cdot)
]_{t_{n}^{i} -}^{t_{n}^{i} +}
\delta(t-t_{n}^{i})
\\
&=& -\sum_{i=0}^\infty[w_n(x,t_n^i)-p_n(t_n^i)] \delta(t-t_n^i)
\chi_{ (-\infty,b(t_n^i)]}(x) \leq0 \qquad\mbox{in }\BbbR\times
(0,\infty).
\end{eqnarray*}
Sending $n\to\infty$ we find that $\cL w\leq0$ and $\cL_1 U\leq0$
in $\BbbR\times(0,\infty)$ in the sense of distributions.
This proves (\ref{3a.eqUw}) and also completes the proof of the proposition.
\end{pf}

\section{Viscosity solutions and boundary crossing
probabilities}\label{sec4}
In this section, we show that the survival distribution $w$ defined
in (\ref{3a.pw})
is the unique viscosity solution of
the
time dependent Kolmogorov forward equation (\ref{WFBP}).
As mentioned above, the
use of viscosity solutions is necessitated by the fact that the
boundaries arising from the
solution to the variational inequality for the inverse boundary
crossing problem do not
have sufficient regularity for us to employ classical solutions. We do
note however, that
when a classical solution exists, it gives the unique viscosity solution.
Consequently, the classical
solution of the partial differential equation (\ref{WFBP}),
if it exists,
is the survival distribution function associated with $b$, defined in
(\ref{3a.pw}).
Since classical solutions of (\ref{UFBP})
are obtained from classical solutions of
(\ref{WFBP}) via the transformation $U=-\partial w/\partial x$, we
also see that a classical solution of (\ref{UFBP}), if it exists, is
the survival probability
density of the first boundary crossing problem that we want to calculate.

For simplicity, we work with $b$ in the class $B_0$ so that the
survival probability associated with $b$ is continuous on
$[0,\infty)$. Furthermore, we work with the function $w$ which is
monotone in the
spatial variable, and smoother than $U$.
The following definition is based on the differential
inequalities/equalities in (\ref{3a.eqUw}).
\begin{definition}\label{def2} Let $b\in B_0$. A \textit{viscosity
solution} (for the survival distribution) of
the boundary crossing problem associated with $b$ is a function $w$ defined
on
$\BbbR\times(0,\infty)$ that has the following properties:
\begin{enumerate}
\item$w\in C(\BbbR\times(0,\infty))$;
$\lim_{t\searrow0} \|w(\cdot,t)+p_0(\cdot,t)-1\|_{L^\infty(\BbbR
)}=0$; $0\leq w\leq1;$

\item$ w(x,t)=w(b(t),t)\ \forall x\leq b(t),t>0;$
$w(x,t)<w(b(t),t) \ \forall x> b(t),t>0$; $\cL w=0$ in $Q_b$;

\item If for a smooth $\varphi$, point $x\in\BbbR$ and time
$t>\delta>0$, the function $\varphi-w$ attains a local minimum at
$(x,t)$ on
$[x-\delta,x+\delta]\times[t-\delta,t]$, then
$\cL\varphi(x,t)\leq0$.
\end{enumerate}
\end{definition}

We define one-sided time derivatives by
\begin{eqnarray*}
\frac{\partial^+\phi(x,t)}{\partial t}&=&\lim_{\Delta t\searrow0}
\frac{\phi(x,t+\Delta t)-\phi(x,t)}{\Delta t},\\
\frac{\partial^-\phi(x,t)}{\partial t}&=&\lim_{\Delta t\searrow
0} \frac{\phi(x,t)-\phi(x,t-\Delta t)}{\Delta t}.
\end{eqnarray*}
Denote by $\cL^\pm$ and $\cL_1^\pm$
the operator $\cL$ and $\cL_1$ with time derivative replaced by the
above one-sided derivative. Then
from the expression of $U_n$ in (\ref{2.Un}), we see that $\cL_1^-
U_n=\cL^-w_n=0$ in $\BbbR\times(t_n^i,t_n^{i+1}]$.
Thus, we have
%
%
\begin{equation}
\label{4.remark} \cL^- U_n(x,t) =0,\qquad \cL^- w_n(x,t)
=0 \qquad\forall (x,t)\in\BbbR\times(0,\infty).
\end{equation}
%
In the uniqueness proof in the following theorem it is
convenient to work with the special case $\sigma\equiv1$. This can be
done without
loss of generality by considering the transformation
%
%
\begin{eqnarray}\label{1.trans}
{\mathbf Y}(x,t)&:=&\int_0^x \frac{1}{\sigma(z,t)}\,dz\qquad \forall x\in
\BbbR,t\geq0,\nonumber\\[-8pt]\\[-8pt]
Y_t&:=&{\mathbf Y}(X_t,t) \qquad\forall t\geq0.\nonumber
\end{eqnarray}
The change of variables $(x,t)\to(y,t)$ via $y={\mathbf Y}(x,t)$ is
smooth and invertible. Also,
by It\^o's lemma,
\[
dY_t 
= \tilde\mu(Y_t,t) \,dt + dB_t.
\]
Here, denoting by $x={\mathbf X}(y,t)$ the inverse of $y={\mathbf Y}(x,t)$,
\[
\tilde\mu(y,t):= -\int_0^x \frac{\sigma_{t}(z,t)}{\sigma
^2(z,t)} \,dz
+\frac{\mu(x,t)}{\sigma(x,t)} - \frac{1}{2}\sigma_{x}(x,t)
\bigg|_{x={\mathbf X}(y,t)}.
\]
Under the transformation, a boundary $b$ for $\{X_t\}$ is transformed
to a
boundary $\tilde b\dvtx t\in[0,\infty)\to{\mathbf Y}(b(t),t)$ for $\{
Y_t\}
$. Similarly, a boundary
$\tilde b$ for $\{Y_t\}$ is transferred back to $b\dvtx t\in[0,\infty
)\to{\mathbf X}(\tilde b(t),t)$.
\begin{theorem}\label{th2.1} Assume that $b\in B_0$.
\begin{enumerate}
\item The survival distribution associated with $b$ defined in
(\ref{3a.pw}) is the unique viscosity solution for the survival
distribution of the
boundary crossing problem associated with $b$.
Consequently, the survival probability $p=\cP[b]$ can be evaluated by
$p(\cdot)=w(-\infty,\cdot)$ where $w$ is the viscosity solution of
the boundary crossing problem associated with $b$ in Definition \ref{def2}.

\item
If $w$ is a classical \textup{(i.e., $w+p_0\in C(\BbbR\times[0,\infty))$,
$\partial_x w\in C(\BbbR\times(0,\infty))$, $\partial_t w,\partial
_x^2w\in C(Q_b)$)}
solution of (\ref{WFBP}),
then $w$ is the survival distribution of the boundary crossing
problem associated with $b$. If $U$ is a classical \textup{(i.e.,
$U-\rho
_0\in C(\BbbR\times[0,\infty))$, $\partial_t U,\partial_x^2 U\in
C(Q_b)$)} solution of (\ref{UFBP}),
then $w(x,t):=\int_x^\infty U(y,t)\,dt$ is
the survival distribution of the boundary crossing problem.
\end{enumerate}
\end{theorem}
\begin{pf}
\textit{Existence}. Let $b\in B_0$ and $\tau, p,w$ be defined as in
(\ref{tauDefinition})--(\ref{3a.pw}).
Then $p=\cP[b]\in C([0,\infty))$. We show that $w$ is a viscosity
solution in the sense of Definition \ref{def2}.
First, we show that $w$ is continuous.
Fix $x\in\BbbR$. For any $t>s>0$,
\begin{eqnarray*}
&&
w(x,t)-w(x,s) \\
&&\qquad=\BbbP(\tau\geq t,X_t>x)-\BbbP(\tau
\geq s,X_s>x)
\\
&&\qquad= \BbbP(\tau\geq t, X_t>x)-\BbbP(\tau\geq t,X_s>x) - \BbbP
(t>\tau\geq s,X_s>x)
\\
&&\qquad= \BbbP(\tau\geq t,X_t>x\geq X_s) - \BbbP(\tau\geq t,
X_s>x\geq X_t)
- \BbbP(t>\tau\geq s,X_s>x).
\end{eqnarray*}
Note that $\BbbP(t>\tau\geq s,X_s>x)\leq\BbbP(t>\tau\geq
s)=p(s)-p(t)$ so we have
\[
|w(x,t)-w(x,s)| \leq\BbbP(X_s>x\geq X_t) +\BbbP
(X_t>x\geq X_s)
+ |p(s)-p(t)|.
\]
Since $p$ is continuous, sending $t\to s$ or $s\to t$ we conclude that
$w(x,\cdot)$ is continuous in $(0,\infty)$. Next, for $x<y$ and $t>0$,
\begin{eqnarray*}
0&\leq& w(x,t)-w(y,t) =\BbbP(\tau\geq t, y\geq X_t>x) \\
&\leq&
\BbbP( y\geq X_t> x)= p_0(y,t)-p_0(x,t).
\end{eqnarray*}
Thus, $w(\cdot,t)$ is continuous, uniform in $t\in[\vep,\infty)$
for any $\vep>0$.
In conclusion, $w\in C(\BbbR\times(0,\infty))$.
Recall from (\ref{3a.ini}) that $\|w(\cdot,t)+p_0(\cdot,t)-1\|
_{L^\infty(\BbbR)}\leq1-p(t)$.
The continuity of $p$ on $[0,\infty)$ then implies that
${\lim_{t\searrow0}} \|w(\cdot,t)+p_0(\cdot,t)-1\|_{L^\infty(\BbbR
)}=\lim_{t\searrow0}(1-p(t))=0$.
Thus, $w$ satisfies the first
requirement of being a viscosity solution.
\begin{remark}
The continuity of the survival probability $p$ plays a central role in
the proof here.
In a more general case, that is, $b\notin B_0$, $w$ is not continuous
so the definition of a
viscosity solution needs to be revised. To avoid such technicalities,
we take the simple case that $b\in B_0$. The work of \citet{CCCS} does
allow discontinuous survival probabilities.
\end{remark}

Note that $\tau(\omega)\geq t$ implies $X_t\geq b(t)$. Hence,
$w(x,t)=\BbbP(\tau\geq t, X_t>x)=w(b(t),t)$ when $x<b(t)$.
Also, since $U=-\partial w/\partial x\geq0$, $U\not\equiv0$ and
$\cL_1 U=0$ in~$Q_b$, we
have $U > 0$ in $Q_b$.
In particular, if $U(x,t)=0$, with $(x,t)\in Q_b$, then the strong
maximum principle
[\citet{FriedmanParabolic}, Theorem 3.5, page 39] implies that
$U(y,t)=0$ for all $y$ such
that $(y,t) \in Q_b$
[and therefore all $y\in\mathbb{R}$, as it is easy to see that
$U(y,t) = 0$ for
$(y,t)\notin Q_b$]. This
is a contradiction, since the Dominated Convergence theorem implies that
\[
0<p(t) \leq p_{n}(t) = \lim_{n\to\infty} \int_{-\infty}^{\infty}
U_{n}(y,t) \,dy
= \int_{-\infty}^{\infty} U(y,t) \,dy
\]
with the application of Dominated Convergence justified by the
bounds $\rho_{0} \geq U_{n} \geq0$ from part 2 of Proposition \ref{prop5}.
Thus, $w(\cdot,t)$ is strictly decreasing in $(b(t),\infty)$, so
$w(x,t)<w(b(t),t)$ for all $x>b(t)$.
Finally, from (\ref{3a.eqUw}), we know $\cL w=0$ in $Q_b$.
Thus, $w$ satisfies the second requirement of being a viscosity solution.
\begin{remark}
Recall that $\cP[b]=\cP[b^*]$ for any boundary $b$;
that is, under our nonconventional definition of default time and
survival probability, both the original boundary $b$ and its
upper-semi-continuous envelope $b^*$ produce the same crossing time,
survival probability,
and survival distribution. Here, we needed to use the
upper-semi-continuous representation of the barrier so $Q_b$ is open
and $w(\cdot,t)$ is strictly decreasing for $x>b(t)$.
\end{remark}

We now verify the third requirement for $w$ being a viscosity solution.
Assume that
$\varphi-w$ attains a local minimum at $(x,t)$ on
$A:=[x-\delta,x+\delta]\times[t-\delta,t]$ where $\varphi$ is smooth
and $t>\delta>0$. We want to show that $\cL\varphi(x,t)\leq0$.
We follow a standard technique for viscosity
solutions. First, we modify $\varphi$ to a new smooth $\psi$ so $\psi
-w$ attains a strict local minimum value zero at $(x,t)$ on $A$.
The function is defined by
\[
\psi(y,s):=\varphi(y,s)+(x-y)^4\delta^{-4}+
(t-s)^2 \delta^{-2 }+ [w(x,t)-\varphi(x,t)].
\]
Then $\psi(x,t)-w(x,t)=0$ and
$\cL\psi(x,t)=\cL\varphi(x,t)$. That $\varphi-w$ attains a local
minimum at $(x,t)$ implies
\begin{eqnarray*}
&&\psi(y,s) -w(y,s)
\\
&&\qquad= (x-y)^4\delta^{-4}+ (t-s)^2
\delta^{-2}
+ [\varphi(y,s)-w(y,s)]-[\varphi(x,t)-w(x,t)] \\
&&\qquad\geq
(x-y)^4\delta^{-4}+ (t-s)^2 \delta^{-2} \qquad\forall
(y,s)\in A.
\end{eqnarray*}
Thus,
$\psi(y,s)-w(y,s)$ attains on $A$ a strict local minimum, being zero, at
$(x,t)$.

Using a standard viscosity solution
technique, the differential inequality $\cL\varphi(x,t)\leq0$ is obtained
by comparison of $\varphi$ with smooth approximations of viscosity
solution candidates. Here, we choose the smooth approximations to be $\{
w_n\}$ introduced in Proposition \ref{prop5}.
For each positive integer $n$, let $w_n$ be defined as in
(\ref{2.wn}) in Proposition \ref{prop5}. Then
$w_n$ is upper-semi-continuous on $\BbbR\times[0,\infty)$, so $\psi
-w_n$ attains a local minimum,
on the closed set $A=[x-\delta,x+\delta]\times[t-\delta,t]$.
On the parabolic boundary of $A$, we have $\psi-w\geq1$ and
$w-w_n>-1$ so $\psi-w_n>0$. At $(x,t)$, $\psi-w_n=w(x,t)-w_n(x,t)\leq
0$. Hence, the minimum is attained in $(x-\delta,x+\delta)\times
(t-\delta,t]$.
We denote by $(x_n,t_n)$ an arbitrary local minimizer of $\psi-w_n$ in $A$.


That $\psi-w_n$ attains a local minimum at $(x_n,t_n)$
implies that at $(x_n,t_n)$, $\partial\psi/\partial x=\partial
w_n/\partial x$,
$\partial\psi/\partial t\leq\partial^- w_n/\partial t$ and
$\partial^2\psi/\partial x^2\geq\partial^2 w_n/\partial x^2$.
Hence,\break $ \cL\psi(x_n$, $t_n)\leq\cL^- w_n(x_n,t_n)=0$.

In order to take the limit, we want to show that $(x_n,t_n)\to(x,t)$
as $n\to\infty$. Intuitively this is obvious since $\psi-w$ attains
a strict local minimum at $(x,t)$ and $w_n\to w$ (uniformly).
Since $\psi(y,s)-w(y,s)\geq(x-y)^4\delta^{-4}+(t-s)^2\delta^{-2}$
with $(y,s)=(x_n,t_n)$,
we have
\begin{eqnarray*}
&&\varlimsup_{n\to\infty}\{
(x_n-x)^4\delta^{-4}+ (t_n-t)^2 \delta^{-2}\} \\
&&\qquad\leq\varlimsup
_{n\to\infty}\{\psi(x_n,t_n)-w(x_n,t_n)\}
\\
&&\qquad=
\varlimsup_{n\to\infty} \{ [\psi(x_n,t_n)-w_n(x_n,t_n)]+
[w_n(x_n,t_n)-w(x_n,t_n)]\}
\\
&&\qquad\leq \varlimsup_{n\to\infty}\Bigl\{
[\psi(x,t)-w_n(x,t)]+\max_{[x-\delta,x+\delta]\times[t-\delta
,t]}|w_n-w|\Bigr\}\\
&&\qquad=\psi(x,t)-w(x,t)=0.
\end{eqnarray*}
Here, we have used the uniform convergence of $w_n\to w$, derived as follows:
\begin{eqnarray*}
0&\leq& w_n(x,s)-w(x,s) \\
&=& \BbbP(\tau_n\geq
s,X_s>x)-\BbbP(\tau
\geq s,X_s>x)
\\
&=& \BbbP(\tau_n\geq s>\tau,X_s>x) \leq \BbbP(\tau_n\geq
s>\tau)
\\
&=& \BbbP(\tau_n\geq s)-\BbbP(\tau\geq s)=p_n(s)-p(s).
\end{eqnarray*}
Thus, we have $\|w_n(\cdot,s)-w(\cdot,s)\|_{L^\infty(\BbbR)}\leq
p_n(s)-p(s)$.
Since $p_n,p$ are continuous, and $p_{n} \searrow p$, the point-wise
convergence of $p_n\to p$ implies local uniform convergence, that is,
$\lim_{n\to\infty}\|p_n-p\|_{L^\infty([0,T])}=0$. Thus,
$\lim_{n\to\infty}\|w_n-w\|_{L^\infty(\BbbR\times[0,T])}=0$ for
any $T>0$.
Hence,
$\lim_{n\to\infty}(x_n,t_n)=(x,t)$. Finally, this implies
$\cL\varphi(x,t)=\cL\psi(x,t)=\lim_{n\to\infty}\cL\psi
(x_n,t_n)\leq
0$.\vspace*{12pt}


\textit{Uniqueness.}
We can assume without loss of generality that $\sigma\equiv1$, since
otherwise we can work with the process $\{Y_t\}_{t\geq0}$ defined in
(\ref{1.trans}). In terms of our viscosity solution, it means that we
make a smooth change of variable $(x,t)\to(y,t)$ via
\[
y=Y(x,t):=\int_0^x\frac{1}{\sigma(z,t)} \,dz.
\]
In the new variables, we are working on the function $w(X(y,t),t)$ and the
barrier is $b(X(y,t))$ where $x=X(y,t)$ is the inverse of $y=Y(x,t)$.
Retaining the notation $(x,t)$ as independent variables,
we can assume that $\cL=\partial_t-\frac12\partial^2_x + \mu
(x,t)\partial_x$.
We denote
\begin{eqnarray*}
M&=&\|\mu\|_{L^\infty(\BbbR\times[0,\infty))}+ \|\partial_x \mu\|
_{L^\infty(\BbbR\times[0,\infty))},\\
R(t)&=&\|\rho_0\|
_{L^\infty(\BbbR\times[t,\infty))}\qquad\forall t>0.
\end{eqnarray*}

We note that $R(t)<\infty$ by the standard Gaussian upper bound on the
fundamental solution
of $\cL_{1}$ [see \citet{FriedmanParabolic}, page 24].
Let $w$ be the
survival probability of the boundary crossing problem. Then $|\partial
_x w|\leq\rho_0$ is uniformly bounded in $\BbbR\times[t_0,\infty)$
for any $t_0>0$. Let $\tilde w$ be an arbitrary viscosity solution of
the boundary crossing problem. We want to show that $w=\tilde w$.

Suppose $w\not=\tilde w$. Then there exists $x_0\in\BbbR,t_0>\delta
>0$ such that either $w(x_0,t_0)>\tilde w(x_0,t_0)+6\delta$ or $\tilde
w(x_0,t_0)>w(x_0,t_0)+6\delta$. In the former case, we set
$(w_1,w_2)=(w,\tilde w)$ and in the latter case we set
$(w_1,w_2)=(\tilde w,w)$. Then both $w_1$ and $w_2$ are viscosity
solutions and
\[
w_1(x_0,t_0)-w_2(x_0,t_0) > 6\delta>0.
\]
By spatial translation, we can assume, without loss of generality,
that
$b(t)<0$ for all
$t\in[0,t_0]$.

We now fix a constant $\vep$ satisfying
\[
0<\vep,\qquad \vep t_0+\vep^4 x_0^4\leq
\min(\delta,1),\qquad \vep^4 +4 \vep^2 M \leq\vep/4.
\]

We need another small positive constant $\eta$ determined as follows.
By the second property of viscosity solutions, we can find $t_1\in
(0,t_0)$ such that
$ \|w_i(\cdot,t_1)+p_0(\cdot,t_1)-1\|_{L^\infty(\BbbR)}<\delta$, $i=1,2$.
Since $p_0(\cdot,t_1)$ is uniformly continuous on $\BbbR$ and $w_2$
is continuous at $(x_0,t_0)$, there exist $\eta_0>0$ such that for
every $\eta\in(0,\eta_0]$, $|w_2(x_0,t_0)-w_2(x_0+\eta,t_0)|\leq
\delta$ and $\|p_0(\cdot,t_1)-p_0(\cdot+\eta,t_1)\|_{L^\infty
(\BbbR)}\leq\delta$. The latter inequality implies
\begin{eqnarray*}
&& \|w_1(\cdot,t_1)-w_2(\cdot+\eta,t_1)\|_{L^\infty
(\BbbR)}\\
&&\qquad\leq\|w_1(\cdot,t_1)+p_0(\cdot,t_1)-1\|_{L^\infty(\BbbR)}
\\
&&\qquad\quad{} + \|w_2(\cdot+\eta,t_1)+p_0(\cdot+\eta,t_1)-1\|
_{L^\infty(\BbbR)}\\
&&\qquad\quad{} +\|p_0(\cdot,t_1)-p_0(\cdot+\eta,t_1)\|_{L^\infty(\BbbR)} \leq3
\delta.
\end{eqnarray*}

Now, we fix an $\eta\in(0,\eta_0]$ such that
\[
0< M \eta[ R(t_1)+ 4\vep^2 ] \leq\vep/4.
\]
%

Consider the continuous function
\[
\Phi(x,t)= w_1(x,t)-w_2(x+\eta,t) -\vep t -\vep^4
x^2,\qquad
x\in\BbbR, t\in[t_1,t_0].
\]
Note that $\Phi(x_0,t_0) =
[w_1(x_0,t_0)-w_2(x_0,t_0)]+[w_2(x_0,t_0)-w_2(x_0+\eta,t_0)]-\vep
t_0-\vep^4 x_0^2 \geq6\delta-\delta-\delta=4\delta$. On the other
hand, when $t=t_1$,
$\Phi(x,t_1) \leq\|w_1(\cdot,t_1)-w_2(\cdot+\eta,t_2)\|_{L^\infty
(\BbbR)} \leq3\delta$.
Hence, there exists
$(x_*,t_*)\in\BbbR\times(t_1,t_0]$ such that $\Phi$ attains at
$(x_*,t_*)$ the global
positive maximum of $\Phi$ on $\BbbR\times[t_1,t_0]$:
\[
\Phi(x_*,t_*)\geq\Phi(x,t)\qquad \forall x\in
\BbbR,
t\in[t_1,t_0].
\]

We consider two cases: (i) $x_*\leq b(t_*)-\eta$; (ii)
$x_*>b(t_*)-\eta$.

Suppose (i) $x_*\leq b(t_*)-\eta$. Then we have
$w_1(x_*,t_*)=w_1(b(t_*),t_*)$ and $w_2(x_*+\eta,t_*)=w_2(b(t_*),t_*)$.
Consequently, since $\|w_2(\cdot,t_*)\|_{L^\infty(\BbbR)}=
w_2(b(t_*)$, $t_*)=w_2(x^*+\eta,t_*)$ and
$x_*< b(t_*)<0$, we obtain
\begin{eqnarray*}
&&\Phi(b(t_*),t_*)-\Phi(x_*,t_*) \\
&&\qquad=\vep^4 \bigl( x_*^2
-b(t_*)^2\bigr)
+ \bigl(w_2(x^*+\eta,t_*)- w_2\bigl(b(t_*)+\eta,t_*\bigr)\bigr) >0
\end{eqnarray*}
contradicting the maximality of $\Phi(x_{*},t_{*})$.


(ii) $x_*>b(t_*)-\eta$. Set
$\varphi(x,t)=w_2(x+\eta,t)+\vep t+\vep^4 x^2$. Then
$\varphi-w_1=-\Phi$ attains at $(x_*,t_*)$ a minimum over $\BbbR
\times[t_0,t_1]$.
Since $x_*>b(t_*)-\eta$, we see that $w_2(\cdot+\eta,\cdot)$ is
smooth in
a neighborhood of $(x_*,t_*)$. Then $\varphi$ is smooth in a small
neighborhood of $(x_*,t_*)$
and $\varphi-w_1$ attains a local minimum at $(x_*,t_*)$.
Since $w_1$ is a viscosity solution, we must have
$\cL\varphi(x_*,t_*)\leq0$.

Now, we calculate $\cL\varphi(x_*,t_*)$. Using $\cL w_2=0$ in $Q_b$
and the fact that $\sigma$ is assumed to be a constant, we have
\begin{eqnarray*}
\cL\varphi(x_*,t_*) &=&\vep- \vep^4 +2\vep^4 x_*\mu
(x_*,t_*) \\
&&{} + [\mu(x_*,t_*)-\mu(x_*+\eta,t_*)] \partial_x w_2(x_*+\eta
,t_*).
\end{eqnarray*}

First, we note that $\Phi(x_*,t_*)\geq\Phi(x_0,t_0)$, so $\vep
t_*+\vep^4 x_*^2 \leq2+\vep t_0+\vep^4 x_0^4\leq3$. This implies
that $\vep^2 |x_*| \leq2$. Hence,
\[
\cL\varphi(x_*,t_*) \geq\vep- \vep^4 - 4\vep^2 M -
M \eta | \partial_x w_2(x_*+\eta,t_*)|.
\]
To estimate $\partial_x w_2(x_*+\eta,t_*)$, we consider two situations.

(a) $w_2=w$ is the survival distribution function. Then $\|\partial_x
w(x_*,t_*)\|\leq\rho_0(x_*$, $t_*)\leq R(t_1)$.

(b) $w_2=\tilde w$. Then $w_2$ is differentiable at $(x_*,t_*)$ and
$w_1(\cdot,t_*)$ is Lipschitz continuous with Lipschitz constant $\|
U(\cdot,t_*)\|_{L^\infty(\BbbR)}$. Hence, sending $h\searrow0$ in
$[(\Phi(x_*\pm h,t_*)-\Phi(x_*,t_*)]/h\geq0$
we derive
\[
|\partial_x w_2(x_*,t_*)|
\leq\|\partial_x w(\cdot,t_*)\|_{L^\infty(\BbbR)} + 2\vep^4 |x^*|
\leq R(t_1)+4\vep^2.
\]

Thus, in either case, we have
\[
\cL\varphi(x_*,t_*) \geq\vep-\vep^4 -4\vep^2 M -
M \eta
\{ R(t_1)+4\vep^2\} \geq\vep/4>0
\]
by our careful choices of $\vep$ and $\eta$. This contradicts $\cL
\varphi(x_*,t_*)\leq0$.
The contradiction implies that we must have $w_1\equiv w_2$. Thus, the
viscosity solution of the boundary crossing problem is unique.\vspace*{12pt}

\textit{Proof of the second assertion.}
The equivalence of classical
solutions of
(\ref{WFBP})
and~(\ref{UFBP})
via $U=-\partial w/\partial x$ is trivial. Here we show that a
classical solution of
(\ref{WFBP})
is a viscosity solution.

Assume that $w$ is a classical solution of (\ref{WFBP}).
Then $U=-\partial w/\partial x$ is a
classical solution of (\ref{UFBP}).
Applying the maximum principle to $U$ and $\rho_0-U$ on~$\bar Q_b$, we
find that $0\leq U\leq\rho_0$. Also, the strong maximum principle
[\citet{FriedmanParabolic}, Theorem
3.5, page 39] shows $U>0$ in $Q_b$ [if $U(x,t)=0$ for $(x,t)\in Q_b$, then
$U(y,s)=0$ for all $(y,s)$ in $Q_b$ with $s\leq t$, contradicting the
initial condition at time
0].
Hence, $w$ is monotone in $x$ and $w(x,t)<w(b(t),t)$ for all $x>b(t),t>0$.
In addition, for each $t_0>0$, comparing $U$ with the solution of $\cL
V=0$ in $\BbbR\times(t_0,\infty)$ with initial value $V(\cdot
,t_0)=U(\cdot,t_0)$ we see that $U\leq V$ on $\BbbR\times[t_0,\infty
)$ so $\int_\BbbR U(y,t) \,dy \leq\int_\BbbR V(y,t) \,dy =\int
_\BbbR U(y,t_0) \,dy$ for every $t>t_0$. This implies that $p(t):=
w(-\infty,t)= w(b(t),t)$ is a decreasing function of $t$.

Next, since $w+p_0\in C(\BbbR\times[0,\infty))$, $\partial_x w\leq
0$, $\partial_x p_0\leq0$, and
$w(\infty,0)=0$ and $w(-\infty,0)=1$, one can show that $\lim
_{t\searrow0}\|w(\cdot,t)+p_0(\cdot,t)-1\|_{L^\infty(\BbbR)}=0$.
Thus, $w$ satisfies the first and second requirements of a viscosity
solution in Definition \ref{def2}.

To verify the third requirement in Definition \ref{def2}, suppose
$\varphi$ is smooth, $x\in\BbbR, t>\delta>0$ and $\varphi-w$
attains a local minimum at $(x,t)$ on $[x-\delta,x+\delta]\times
[t-\delta,t]$. We want to show that $\cL\varphi(x,t)\leq0$. We
consider two cases:
(i) $x>b(t)$ and (ii) $x\leq b(t)$.

\mbox{}\hphantom{i}(i) Suppose $x>b(t)$. Then $(x,t)$ is an interior point of $Q_b$,
in which $w$ is smooth. Since $\varphi-w$ attains a local minimum at
$(x,t)$ on $[x-\delta,x+\delta]\times[t-\delta,t]$, we have
$\partial_t\varphi(x,t) \leq\partial_t w(x,t), \partial_x \varphi
(x,t)=\partial_x w(x,t)$ and
$\partial^2_x \varphi(x,t)\geq\partial^2_x w(x,t)$. This implies
that $\cL\varphi(x,t)\leq\cL w(x,t)=0$.

(ii) Suppose $x\leq b(t)$. Then $\partial_x \varphi(x,t)=\partial_x
w(x,t)=0$. Note that $\varphi(x-z,t)-\varphi(x,t)\geq
w(x-z,t)-w(x,t)=0$ for all
$z\geq0$. This implies, since $\varphi$ is a smooth and $\partial
_x\varphi(x,t)=0$, that $\partial^2_{xx}\varphi(x,t)\geq0$.

To complete the proof that $\cL\varphi(x,t)\leq0$, it suffices to
show that $\partial_t\varphi\leq0$. Suppose on the contrary that
$\partial_t\varphi(x,t)>0$. Then there exists $\vep\in(0,\delta)$
such that
$\varphi(x,t-s)<\varphi(x,t)$ for all $s\in(0,\vep]$. As $\varphi
-w$ attains a local minimum at $(x,t)$, we see that $w(x,t-s)\leq
w(x,t)-\varphi(x,t)+\varphi(x,t-s)<w(x,t)=p(t)\leq p(t-s)$ for all
$s\in(0,\vep]$. Thus, $[x,\infty)\times[t-\vep,t)\subset Q_b$.
Since $w$ is monotone, we also have $w(y,t-s)\leq p(t)=w(x,t)$ for all
$y>x,s\in[0,\vep]$. That is,
$w$ attains at $(x,t)$ a local maximum over the region $[x,\infty
)\times[t-\vep,t]$.
Hence, applying Hopf's lemma [\citet{ProtterWeinberger}, Theorem
3.3]
for $w$ on $[x,\infty)\times(t-s,t)$, we have $w_x(x,t-)<0$, which
contradicts the definition of a classical solution that $\partial_x
w(x,t-)=\partial_x w(x,t)=0$.
Thus, we must have $\partial_t \varphi(x,t)\leq0$.
Together with $\partial_x \varphi(x,t)=0,\partial^2_{xx}\varphi
(x,t)\geq0$, we conclude that $\cL\varphi(x,t)\leq0$.

Hence, $w$ is a viscosity solution.
Applying the conclusion of the first assertion, we then see that $w$ is
the survival
distribution of the boundary crossing problem associated with $b$.
\end{pf}
\begin{remark}\label{re4.1} 
(i) The addition of the term $\vep^4 x^2$ confines our attention in searching
for the maximum of $\Phi$ to a compact set $[-2\vep^{-2},2\vep
^{-2}]\times[0,t_0]$.
So we indeed only need $w+p_0\in C(\BbbR\times[0,\infty))$ and
$w+p_0-1=0$ on $\BbbR\times\{0\}$ in the definition of viscosity solutions.

(ii)
Fix $t_2>t_1\geq0$. In the proof, if we take $w_2$ to be a solution of
$\cL w_2=0$ in $\BbbR\times(t_1,t_2)$ subject to initial condition
$w_2(\cdot,t_1)\geq w(\cdot,t_1)$. Then following the proof we see
that $\sup_{\BbbR\times[t_1,t_2)}(w-w_2)>0$ is impossible. Thus, we
have $w\leq w_2$ on $\BbbR\times[t_1,t_2)$. This is a simple version
of the comparison principle in
the theory of viscosity solutions. This result will be used in the next section.
\end{remark}

\section{Viscosity solutions for the inverse boundary crossing problem}

The inverse boundary crossing problem is to find $b$, for a given
$p$, such that $p$ is the survival probability associated with $b$. In
this section,
we prove that for any $p\in P_{0}$,
from the viscosity solution of the variational inequality (\ref{2.var}),
we can find an unique $b\in B_{0}$ such that the resulting
$p$ gives the survival distribution of the
first time that $X$ crosses $b$.
Since the forward problem maps $B_0$ to $P_0$, we study,
for simplicity, the inverse problem for $p\in P_0$, though in \citet
{CCCS} the
variational inequality (\ref{2.var}) was considered for
more general survival functions.

\subsection{Viscosity solutions}
In general, classical solutions of the variational inequality
(\ref{2.var}) for the inverse problem may not exist.
In \citet{CCCS}, viscosity solutions were introduced, and it was
shown that for any $p$ satisfying
%
%
\begin{equation}\label{2.CCCS}
p(0+)=1\geq p(s)=p(s-)\geq p(t)\geq0 \qquad\forall
t>s>0,
\end{equation}
there exists a unique viscosity
solution. From this solution, we can define a boundary
$b$ such that $Q_b=\{w<p\}$, and consider it as a candidate for the
solution to the inverse
boundary crossing problem. To verify that $b$ is indeed a solution, we
show that $w$ is a
viscosity solution to the direct problem (\ref{WFBP}), and then appeal to
Theorem \ref{th2.1} to see that $w$ and $p$ give the survival
distribution of the
first time that $X$ crosses $b$.

When $p\in P_0$, we know a priori that the unique solution of
the variational inequality is continuous so many technicalities in
\citet{CCCS}
regarding the definition, existence, and uniqueness of viscosity
solutions can be ignored.
In particular, the viscosity solution introduced in \citet{CCCS} can be
reformulated (removing those specifics that take care of
discontinuities) as follows.
\begin{definition}\label{def1}
Let $p\in P_0$ be given. A \textit{viscosity solution for the survival
distribution of the
inverse boundary crossing problem} associated with
$p$ is a function $w$ defined on $\BbbR\times(0,\infty)$ that has
the following properties:
\begin{enumerate}
\item$w+p_0\in C(\BbbR\times(0,\infty))$, $\lim_{t\searrow0}\|
p_0(\cdot,t)+w(\cdot,t)-1\|_{L^\infty(\BbbR)}=0$;

\item$0\leq w(x,t) \leq p(t)$,
$\forall(x,t)\in\BbbR\times(0,\infty)$ and
$\cL w(x,t) =0
$ in the set
$Q:=\{(x,t) | t>0, w(x,t)<p(t)\}$;

\item if $x\in\BbbR$ and $ t>\delta>0$, and $\varphi$ is a
smooth function such that $\varphi-w$ attains at
$(x,t)$ a local minimum on $[x-\delta,x+\delta]\times[t-\delta,t]$,
then $\cL\varphi(x,t)\leq0$.
\end{enumerate}

A \textit{viscosity solution of the inverse boundary crossing problem}
associated with $p$ is the function $b$ given by
%
%
\begin{equation}\label{2.ba}
b(t):=\inf\{ x\in\BbbR| w(x,t)<p(t)\}\qquad
\forall t>0,
\end{equation}
where
$w$ is a viscosity solution for the survival distribution of the
inverse boundary crossing problem associated with
$p$. If there is a unique viscosity solution, we denote $b=\cB[p]$.
\end{definition}

The remainder of this section is devoted to a proof of the main result
of the paper,
Theorem \ref{MainTheorem}, stated in the \hyperref[intro]{Introduction}.
%
%
\begin{pf*}{Proof of Theorem \ref{MainTheorem}}
The fourth assertion follows from the first assertion and the following
facts which are easy to verify: a classical solution
of (\ref{2.var}) is automatically a viscosity solution, and if $(U,b)$
is a
classical solution of (\ref{UFBP}), then $(w,b)$ with $w$
defined by $w(x,t)=\int_x^\infty U(y,t)\,dy$ is a classical solution of
(\ref{2.var}).
We divide the proof of the first three assertions into several
parts.\vspace*{12pt}

\textit{Existence and uniqueness of a viscosity solution.}
The proof of the existence of a unique
viscosity solution, together with the properties\vspace*{1pt}
(\ref{monotonew}), is the major result of \citet{CCCS} and hence is
omitted here.\footnote{In \citet{CCCS}, it was assumed that $\mathbb
{P}(X_{0}=0)=1$. The
same techniques can be applied to prove existence and uniqueness for
more general
initial distributions considered here.}
It is important to note that, by the monotonicity of $w$ in the spatial
variable and
the definition of $b$ in (\ref{2.ba}), we have
\[
Q:=\{(x,t) | t>0,w(x,t)<p(t)\} = Q_b:=\{(x,t) |
t>0,x>b(t)\}.\vspace*{12pt}
\]


\textit{Weak regularity of the free boundary.}
The regularity of the free boundary $b=\cB[p]$ defined by (\ref
{2.ba}) was
not discussed in \citet{CCCS}. Here, under the assumption that $p\in
P_0$, we establish a very basic regularity result on $b$.
In particular, we show that $b\in B_0$.

We begin by showing that $b(t) < \infty$ for every $t>0$.
Indeed, if $b(t)=\infty$, then by the definition $b(t)=\inf\{x |
w(x,t)<p(t)\}$ we see that $w(x,t)=p(t)$ for all $x\in\BbbR$. Since
$\lim_{x\to\infty} w(x,t)=0$ (recalling $w\leq1-p_0$), we see that
$w(\cdot,t)\equiv0$. This contradicts the assumption that $p\in P_0$,
since $p\in P_0$ guarantees $p(s)>0$ for all $s\in[0,\infty)$.

Next, we show that $X_0\geq b^{*}(0) =\varlimsup_{t\searrow0}b(t)$
almost surely.
To see this, we use the estimate $w\leq1-p_0$ to derive
$p(t)=w(b(t),t)\leq1-p_0(b(t),t)$.
This implies that $\varlimsup_{t\searrow0} p_0(b(t),t)\leq\lim
_{t\to0} (1-p(t))=0$ since
$p\in P_0$ gives $p\in C([0,\infty))$ and $p(0)=1$.
Now suppose to the contrary that $\BbbP(X_0<b^{*}(0))>0$. Then there
exists $\delta>0$ and
$\vep>0$ such that
$p_0(b^{*}(0)-2\delta,0)=\BbbP(X_0\leq b^{*}(0)-2\delta)>2\vep$.
Consequently,
there exists $t_0>0$ such that $p_0(b^{*}(0)-\delta,t)=\BbbP
(X_t<b^{*}(0)-\delta)>\vep$ for all $t\in[0,t_0]$.
However, by the definition of $b^{*}(0)$, there exists a sequence
$t_k\searrow0$ such
that $\lim_{k\to\infty}
b(t_k)=b^{*}(0)$. For all sufficiently large $k$, we have
$b(t_k)>b^{*}(0)-\delta$ which implies that
$p_0(b(t_k),t_k)\geq p(b^{*}(0)-\delta,t_k) >\vep$,\vadjust{\goodbreak} so we obtain
$\varlimsup_{t\searrow0} p_0(b(t),t)\geq\vep$. This contradicts
$\lim_{t\searrow0} p_0(b(t),t)=0$.
Hence, we must have $\BbbP(X_0<b^{*}(0))=0$, that is, $X_0\geq
b^{*}(0)$ a.s.

The next step is to show that $b:=\cB[p]$ is upper-semi-continuous on
$[0,\infty)$.
First of all, the definition $b(0):=\limsup_{t\searrow0} b(t)$
implies that $b$ is upper-semi-continuous at $t=0$. Next, let $t>0$ be
arbitrary. We consider two cases: (i) $b(t)>-\infty$, (ii)
$b(t)=-\infty$.

\mbox{}\hphantom{i}(i) Suppose $b(t)>-\infty$.
Fix any $\vep>0$. Then $p(t)-w(b(t)+\vep,t)>0$. By continuity,
$p(s)-w(b(t)+\vep,s)>0$ for all $s$ in a neighborhood of $t$. This
implies that $b(s)<b(t)+\vep$
for all $s$ in a neighborhood of $t$. Consequently,
$\varlimsup_{s\to t}b(s)\leq b(t)+\vep$. As $\vep>0$ is arbitrary,
we have
$\varlimsup_{s\to t} b(s)\leq b(t)$.

(ii) Suppose $b(t)=-\infty$. Then for any $M>0$, we have
$p(t)-w(-M,t)>0$. Consequently,
$p(s)-w(-M,s)>0$ for all $s$ sufficiently close to $t$. Hence,
$b(s)<-M$ for all $s$ sufficiently close to $t$. This implies that
$\varlimsup_{s\to t} b(s) \leq-M$. As $M$ can be made arbitrarily
large, we hence see that
$\lim_{s\to t} b(s) =-\infty=b(t)$.

In conclusion, $b\dvtx[0,\infty)\to[-\infty,\infty)$ is
upper-semi-continuous.

Let $Q:=\{w<p\}:= \{(x,t)\in\BbbR\times(0,\infty) | w(x,t)<p(t)\}
$ and
$Q_b:=\{(x,t)\in\BbbR\times(0,\infty) | x>b(t)\}$. Then $Q=Q_b$
and $U:=-\partial_x w>0$ in $Q_b$.
Indeed, $Q=Q_b$ follows from the definition of $b=\cB[p]$ in (\ref
{2.ba}) and the monotonicity of $w$ in the spatial variable [see \citet
{CCCS}]. In addition, since $b$ is upper-semi-continuous and bounded
above in any compact interval in $[0,\infty)$, any points
$(x_{1},t_{1})$, $(x_{2},t_{2})$,
$t_{1}\leq t_{2}$ in $Q_{b}$ can be connected by a smooth curve
$x=h(t)$ in $Q_b$ defined on $t\in[t_1,t_{2}]$ such that $h(t_i)=x_i$.
Applying the strong maximum principle [\citet{FriedmanParabolic},
Theorem 3.5,
page 39] to $U:=-\partial_x w$
($U\geq0$ by the monotonicity of $w$), we conclude that $U>0$ in
$Q=Q_b$ [an
elementary argument shows there cannot exist a $t_2>0$ such that
$U(\cdot, t_{1}) \equiv0$ in $Q_{b}$ for all $t_{1} \leq t_{2} $].

Next, we show that $b=b^*=b^*_-$. Upper-semi-continuity ($b=b^*$) was
shown above,
so it remains to prove that $b(t)=\varlimsup_{s\nearrow t}
b(s)=:b^*_-(t)$ for every $t>0$.
Let $t>0$, and suppose $b(t)\not= b_-^*(t)$. Then
$b(t)=b^*(t)>b^*_-(t)$. Set
$\delta=[b(t)-b^*_-(t)]/4$. By the definition of $b^*_-(t)$, we can find
$\vep>0$ such that
$b(s)< b^*_-(t)+\delta$ for all $s\in[t-2\vep,t)$.
Then $D:=[b^*_-(t)+\delta,b(t)]\times[t-2\vep,t)$ is a subset of
$Q:=\{w<p\}$.
Since we know that $\cL_1 U=0$ and $U=-\partial w/\partial x>0$ in
$Q$. We can apply the
Harnack inequality on the cube $(b^*_-(t)+\delta,b(t))\times(t-2\vep
,t)$ to conclude that there exists a positive constant $\eta>0$ such
that
$U>\eta$ in $[b^*_-(t)+2\delta, b(t)-\delta]\times(t-\vep,t)$.
Consequently,
\begin{eqnarray*}
&&
w\bigl(b^*_{-}(t)+2\delta,s\bigr)-
w\bigl(b(t)-\delta,s\bigr) \\
&&\qquad=\int_{b^*_{-}(t)+2\delta}^{b(t)-\delta} U(y,s)\,dy
\geq[b(t)-b^*_-(t)-3\delta]\eta\\
&&\qquad= \delta\eta\qquad\forall
s\in(t-\vep,t).
\end{eqnarray*}

Sending $s\nearrow t$ we then conclude that
$w(b^*+2\delta,t)\geq
w(b(t)-\delta,t)+\delta\eta\geq p(t)+\delta\eta$, which
violates the requirement that $w\leq p$ for a viscosity
solution. Hence, we must have
$b(t)=b^*_-(t)$.
In summary, $b=b^*=b^*_-$. This also implies that
\[
b(t)=\varlimsup_{s\to t} b(s)=\varlimsup_{s\nearrow t} b(s)\geq
\varlimsup_{s\searrow t} b(s).
\]

Finally, to show that $b\in B_0$, it remains to show that the survival
probability $\tilde p:=\cP[b]$ associated with $b$ has the property
$\lim_{t\searrow0} \tilde p(t)=1$.
For this, we consider the sequence $\{w_n\}$, associated with $b$,
defined in Proposition \ref{prop5}.
It follows from a (viscosity solution) comparison principle,
applied iteratively to $\BbbR\times(t_n^i,t_n^{i+1}]$ ($t_n^0:=0$)
for $i=0,1,\ldots,$ that $w_n\geq w$; see Remark \ref{re4.1}. Taking
the limit,
we find that $w\leq\lim_{n\to\infty} w_n$. This implies that
$p(t)=w(-\infty,t)\leq\lim_{n\to\infty} w_n(-\infty,t) = \tilde p(t)$.
Since our assumption $p\in P_0$ implies that $\lim_{t\searrow
0}p(t)=1$, we also know that $\lim_{t\searrow0}\tilde p(t)=1$. Thus,
we have shown that $b\in B_0$.\vspace*{12pt}

\textit{Verification that the boundary derived from the
variational inequality has
the required crossing time distribution.}
Given $p\in P_{0}$, let $b=\cB[p]$ be the boundary derived from the
unique solution of
the variational inequality (\ref{2.var}). We need to show that $\cP
[b]=p$, that is, that $b$ is
truly a solution of the inverse boundary crossing problem.
Summarizing, this means that we want to show that
$(\cP\circ\cB)[p]=p$ for every $p\in P_0$.

Let $w$ be the unique \textit{viscosity} solution for survival
distribution of the inverse problem associated with $p$ as given in
Definition \ref{def1}.
Define $b=\cB[p]$ as in
(\ref{2.ba}). Let $\tilde p=\cP[b]$. We want to show that $\tilde
p=p$. It is enough to show that $w$ is a viscosity solution of the survival
probability for the boundary crossing problem associated with $b$,
since in this case
part 1. of Theorem \ref{th2.1} yields that $\tilde p(t) =
w(-\infty,t)$, while taking
limits as $x$ goes to $-\infty$ in (\ref{monotonew}) gives that
$p(t)=w(-\infty,t)$.
By checking the Definitions \ref{def1} and \ref{def2} of viscosity
solutions, one readily sees that $w$ being a viscosity solution in the
sense of Definition \ref{def1} implies that $w$ is indeed the
viscosity solution
in the sense of Definition \ref{def2}, provided that $Q:=\{w<p\}
=Q_b:=\{(x,t) | x>b(t),t>0\}$.
But this last property is immediate from the monotonicity of $w$.
We thus conclude that $w$ is indeed the viscosity solution of the
survival distribution of the
boundary crossing problem
associated with $b$. Consequently,
$\cP[b](t)=w(-\infty,t)=p(t)$, so we have $\cP[b]=p$ and $p=\cP
[b]=(\cP\circ\cB)[p]$.\vspace*{12pt}

\textit{Uniqueness of the solution of the inverse boundary crossing
problem in the class~$B_0$.}
For a given $p\in P_0$, we have shown that $\cB[p]$ is a solution of
the original inverse
boundary crossing problem. Here we show that there is indeed only one
such $b$ in the class
$B_0$. To show this, it suffices to show that
$(\cB\circ\cP)[ b]= b$ for every $b\in B_0$, since this implies that
if $\cP[\tilde b] = p$ then
$\tilde b=(\cB\circ\cP)[\tilde b]=\cB[p]$.

Let $\tilde b\in B_0$. Define $(\tau,p,w)$ as in (\ref
{tauDefinition})--(\ref{3a.pw}).
That is, $p=\cP[\tilde b]$ and $w$ are the survival probability and
distribution of the
boundary crossing problem associated with $\tilde b$.
Since $\tilde b$ is upper-semi-continuous we can derive from
Proposition \ref{newth2} and the strong maximum principle that
$U:=-\partial_x w>0$ in $Q_{\tilde b}$ (see the proof of
Theorem~\ref{th2.1}).
This implies that
$Q_{\tilde b} \subset\{ w<p\}$.
Also, since $w(x,t)=\BbbP(\tau\geq t, X_t>x)$
we know that $w(x,t)=p(t)$ when $x\leq\tilde b(t)$. Thus, $Q_{\tilde
b}=\{w<p\}$.

By Theorem \ref{th2.1}, $w$ is a \textit{viscosity} solution of
the survival probability distribution of the boundary crossing problem
in the
sense of Definition \ref{def2}, associated with~$\tilde b$. By
checking the definition
of a viscosity solution of the variational inequality associated with
$p$ (Definition \ref{def1}),
we find that $w$ is indeed a viscosity solution associated with $p=\cP
[\tilde b]$.
Now, according to the definition of $\cB[p]$ in (\ref{2.ba}),
$\cB[p](t) = \inf\{ x | w(x,t)<p(t)\}$. Since $Q_{\tilde b}=\{w<p\}
$, we see that $\cB[p]=\tilde b$.
Thus $\tilde b=\cB[p]=(\cB\circ\cP)[\tilde b]$ for every $\tilde
b\in B_0$.
This completes the proof of Theorem \ref{MainTheorem}.
\end{pf*}
%

\subsection{Continuity of the free boundary in the inverse boundary
crossing problem}

In this subsection, we investigate the continuity of the free boundary
$b=\cB[p]$
for the inverse boundary crossing problem for
$p\in P_0$. We already know that $b$ is upper-semi-continuous, and
since $b=b^{*}_{-}$,
it cannot ``jump up.'' For $b$ to be continuous, we need to prevent
it from ``jumping down.'' Note the fact that if $p$ is a constant in
an open interval, then $b=-\infty$ in that interval. Hence, to
eliminate steep drops of $b$ we require a lower bound on the rate of
decrease of $p$.
We consider the following:
\[
L(p,T_1,T_2) := \inf_{T_1\leq s<t\leq T_2}\frac
{p(s)-p(t)}{t-s}
\qquad\forall 0\leq T_1<T_2.
\]
The following proposition gives a sufficient condition for the boundary
to be continuous,
in the case that $X$ is a standard Brownian motion.
\begin{proposition} Suppose that
$X$ is a standard Brownian motion, that is, \mbox{$\mu\equiv0$},
$\sigma\equiv1$, and $p_{0}(x,0) = \chi_{[0,\infty)}(x)$.
Let $p\in P_0$ and $b=\cB[p]$.

1. If
$L(p,T_1,T_2)>0$ for some positive $T_1,T_2$ with $T_1<T_2$,
then $b=\cB[p]$ is continuous on
$(T_1,T_2)$.


2. Assume that $L(p,0,T) > 0$ for every $T>0$.
Then $b\in C([0,\infty))$.
\end{proposition}
\begin{pf}
1.
Let $t_1\geq0$ be arbitrary. Define
\begin{eqnarray*}
\tilde w(x,t)&=&\frac{ w(b(t_1)+x,
t+t_1)}{p(t_1)},\qquad \tilde
b(t)= b(t_1+t)-b(t_1), \\
\tilde p(t)&=&\frac{p(t+t_1)}{p(t_1)}.
\end{eqnarray*}
Then $(\tilde w,\tilde b)$ is the solution of the inverse problem with
initial value $\tilde w(\cdot,0)$ and survival probability $\tilde p$. This
statement follows by an immediate application of the definitions. Note that
$\tilde w(x,t) = \mathbb{P}(X_{t+t_{1}} > x, \tau\geq t+t_{1} | \tau
\geq t_{1})$.
$\tilde p(t) = \mathbb{P}(\tau\geq t+t_{1} | \tau\geq t_{1})$. The
conditional probabilities
and the boundary from time $t_{1}$ on are the same as the solution of
the inverse problem started with the initial position equal to the
conditional distribution of
$X$ given that $\tau\geq t_{1}$.

Now let $(\underline{w},\underline{b})$ be the solution of the
inverse problem with initial data
$\chi_{(-\infty,0)}$ and survival probability $ \tilde p$. Note that
$\underline{w}(x,0)=1=\tilde w(x,0)$ for $x<0$ and $\underline
w(x,0)=0\leq w(x,0)$ for $x\geq0$.
Hence, $\underline{w}(\cdot,0)\leq\tilde w(\cdot,0)$. It then
follows from a comparison principle
[cf. the proof of Lemma 4.2 in \citet{CCCS}] that $\underline{w}\leq
\tilde w$ and that
$\underline b\leq\tilde b$. Again, this is obvious from the
probabilistic interpretation of
the problem. The boundary $\underline{b}$ is the one that produces the
hitting distribution
$\tilde p$ when the process starts at $b(t_{1})$ at time $t_{1}$. The
boundary $\tilde b$
produces the same hitting distribution with the process started at the
conditional distribution of
$X_{t_{1}}$ given that $\tau\geq t_{1}$. Since in this case we must
have $X_{t_{1}} \geq
b(t_{1})$, we have that the boundaries $\tilde b$ and $\underline{b}$
produce the same
hitting distribution for the process $X$, with $\tilde b$ arising from
$X$ starting at a higher point
with probability 1. Therefore, we must have $\tilde b \geq\underline{b}$.

Thus, for $0<t\leq1/2$,
\[
b(t_1+t)-b(t_1) =\tilde b(t) \geq\underline{b}(t)
\geq-[1+O(1)] \sqrt{-2t \log[1-\tilde p(t)]}
\]
by the estimate on line 16, page 867 of \citet{CCCS}.
Upon noting that
\begin{eqnarray*}
|{\log}[1-\tilde p(t)]| &=& |{\log}[p(t_1+t)-p(t_1)]-\log p(t_1)|
\\
&=& \bigl|\log\bigl(t\dot p(t_1+\theta) \bigr)\bigr|+O(1)\\
&=& |{\log t}|+O(1) = O(1) |{\log t}|,
\end{eqnarray*}
where
\[
\dot p(t) := \limsup_{s\nearrow t} \frac{p(t)-p(s)}{t-s} \in
[0,\infty]
\]
we find that there exists a constant $C(t_1)$ such that
\[
b(t+t_1)-b(t_1) \geq-C(t_1) \sqrt{t |{\log t}|}\qquad
\forall
t\in(0,1/2].
\]


Now pick any $t\in(T_1,T_2)$. Let $\{t_i\}_{i=1}^\infty$ be a
sequence in $[T_1,t)$ such that $\lim_{i\to\infty} t_i=t$ and
$\lim_{i\to\infty} b(t_i)=b(t)$ [recalling
$b(t)=b^*_-(t):=\varlimsup_{s\nearrow t}b(s)$]. Then setting
$h_i=[t-t_i]/2$, we have
\[
b(t_i+h)>b(t_i)-c \sqrt{|h_i\log h_i|} \qquad\forall
h\in[h_i,3h_i].
\]
This implies that
\[
\inf_{s\in[t-h_i,t+h_i]} b(s) \geq b(t_i)-c\sqrt
{h_i|{\log
h_i}|},
\]
so that
\[
\varliminf_{s\to t} b(s)\geq b(t_i)-c\sqrt{h_i|{\log h_i}|}.
\]
Sending $i\to\infty$ we then obtain $\varliminf_{s\to t}
b(s)\geq b(t)$. Thus, $b$ is lower-semi-continuous in
$(T_1,T_2)$. Since
$b$ is also upper-semi-continuous, we see that
$b$ is continuous in $(T_1,T_2)$.


2. By the first assertion, we know that $b$ is continuous in
$(0,\infty)$.
%
At $t=0$, since $p_0(x,0)>0$ for all $x>0$, the
proof in \citeauthor{CCCS} [(\citeyear{CCCS}), Lemma 4.5, page~865] implies that there exists
a positive
constant $C$ that depends on $L(p,0,1/2)$ such that
\[
b(t) \geq-C \sqrt{t|{\log t}|} \qquad\forall t\in
[0,1/2].
\]
This implies that $\varliminf_{t\searrow0} b(t)\geq0$.

We recall that
$1-w(x,t)-p_0(x,t) \geq0$. Evaluating this inequality at $x=b(t)$ gives
$p_0(b(t),t) \leq1-p(t)$ for all $t>0$.
Since $p_0(x,0)>0$, sending $t\searrow0$ we conclude that
$\varlimsup_{t\searrow0} b(t)\leq0$. Thus,
$b(0):=\lim_{t\searrow0}b(t)=0$. This completes the proof.
\end{pf}

Note that if the hitting
time density $-\dot p$ is everywhere strictly positive, then we obtain
that $b = \cB[p]$ is
continuous. In particular, this criterion is satisfied by the boundary
arising from the exponential
distribution and hence provides the solution to the inverse boundary
crossing problem as originally
proposed by A. N. Shiryayev.

%

%
\printaddresses

\end{document}